\documentclass{article}

\usepackage{PRIMEarxiv}

\usepackage[utf8]{inputenc} 
\usepackage[T1]{fontenc}    
\usepackage{hyperref}       
\usepackage{url}            
\usepackage{booktabs}       
\usepackage{amsfonts}       
\usepackage{amsmath} 
\usepackage{nicefrac}       
\usepackage{microtype}      
\usepackage{lipsum}
\usepackage{fancyhdr}       
\usepackage{graphicx}       
\usepackage{color}
\usepackage{tikz}
\usepackage{multirow}                                   
\usepackage{longtable}                                  
\usepackage{lscape}                                     
\usepackage{xcolor,colortbl} 
\usepackage[linesnumbered,ruled,vlined]{algorithm2e}

\usepackage{url} 
\usepackage{geometry}
\usepackage{graphicx}
\usepackage{subfigure}
\usepackage{pdfpages}
\usepackage{lmodern}
\usepackage{nomencl}
\usepackage{ragged2e}
\usepackage{fancyhdr}
\usepackage{float}
\usepackage{verbatim}
\usepackage{natbib}
\usepackage[skip=5pt]{caption}

\title{A study of column generation embedded in  scalarization methods for the bi-objective cutting stock problem
}

\author{
  Jennifer C. Borges\\
  Escola Municipal Profa. Josiany França\\
  Prefeitura Municipal de Uberlândia\\
  Uberlândia, 38412-402, MG,  Brazil.\\
  jennifercristinamat@gmail.com \\
   \And
   Helenice de O. Florentino\\
   Universidade Estadual Paulista (Unesp)\\
   IBB, Departamento de Bioestatística\\
   Campus  Botucatu, 18618-689, SP, Brazil\\
   helenice.silva@unesp.br\\
   \And
   Socorro Rangel\\
   Universidade Estadual Paulista (Unesp)\\
   IBILCE,  Departmento de Matemática\\
   Campus São José do Rio Preto, 15054-000, SP, Brazil\\
   socorro.rangel@unesp.br (Corresponding author)
}

\begin{document}
\maketitle

\begin{abstract}
Research on multi-objective combinatorial optimization and on the Cutting Stock Problem (CSP) has been widely developed over the years. In contrast, the multi-objective Cutting Stock Problem has received limited attention and has been explored in only a small number of studies. In this paper a bi-objective study of the one-dimensional and the two-dimensional CSP is presented. It distinguishes itself from other research in the literature in two key aspects, among others. The first regards the model used to represent the problem and the second is the solution strategy based on dynamic column generation embedded into scalarization methods. Three methods adapted from the literature to analyse the trade-off between the minimization of the total number of objects and the total number of saw cycles are implemented. The computational results show that the use of dynamic column generation provides a better approximation of the Pareto front.  As for the scalarization methods, none stood out. In fact, they can be viewed as complementary. The cardinality and hypervolume of an approximation of the Pareto front built from the union of the points generated by the three methods  is always greater than the ones generated by each one. Dealing with the  multi-objective nature of CSP, this study provides insights and computational tools to help industry obtain pragmatic solutions.

\end{abstract}


\section{Introduction}\label{sec:intro}
Making decisions while acknowledging that there is not a single solution
that attends a set of objectives is a challenge. The cutting stock
problem (CSP) is an example of such a situation, although it is
frequently considered a mono-objective optimization problem (\textit{e.g.}
\cite{Morabito2009, Oliveiraetal2023, Cherry2023, Kelly2023}), even when two or more objectives are
clearly stated (\textit{e.g.} \cite{Malaguti}). The most frequent
criterion used to solve the CSP is minimization of total
waste, followed by maximization of the productivity of the cutting machine. The latter can be achieved by reducing the number of
different cutting patterns (minimizing setup \emph{\textit{e.g.}}
\cite{araujoPoldi2014} ) or taking into account that some cutting
machines allow the objects to be stacked so that they can be cut
simultaneously (minimizing saw cycles \cite{Yanasse2008}).

The literature on multi-objective combinatorial optimization is quite
extensive, an interesting tutorial on the subject is given in
\cite{Ehrgott2003}. Several approaches have been used to
solve the problem either heuristically (\textit{ e.g.}\cite{Soylu2015}) or exactly
(\textit{e.g.}\cite{Boland2015}). Solving the CSP explicitly addressing
conflicting objectives has received less attention. For the
one-dimensional case, the papers \cite{Angelo2018, Angelo2021, araujoPoldi2014, Kolen} consider the objectives of the minimization
of total waste and setup; in  \cite{Leduino} sufficient conditions are presented  taking into account  the minimization of the total number of {objects} and  setup and in \cite{Rangel2} the two  objectives considered are the minimization of the total number of {objects} and  the total number of cycles. Respicio and Captivo  \cite{Respicio} address the problem of cutting pattern sequencing tanking into account the
minimization of the maximum number of open stacks and the minimization of the average order spread. A trade-off among the minimization of the maximum number of open stacks, total waste and setup is considered in \cite{Munoz}. In \cite{Salles} the last two objectives as well as the minimization of saw cycles are considered. In \cite{Sinuany-Stern-Weine} a bi-ojective 1D-CSP is adressed considering waste minimization as the main optimization criterion and  the  maximization of accumulated  leftovers for use in the future as the secondary objective. The 1D-CSP is also addressed in \cite{wascher1990lp}  exploring the cut of  steel plates rolled up into coils and considering several optimization criteria, among them the minimization of  costs of material, storage, excess of production, left-overs and trim-loss. { The multiperiod bi-objective 1D-CSP considering the minimization of total waste and total costs related to inventory of objects and items is addressed in  \cite{PoldiMO}}.

For the two-dimensional CSP,  the irregular case considering the minimization of the total waste and of the total time to generate cutting patterns is addressed in \cite{Gomez}. The regular case is addressed in  \cite{Toscano, Ahmed}. Toscano \textit{et al } \cite{Toscano} consider the minimization of  the total  number of objects  and of the total number of saw cycles in the context of a furniture industry. Mellouli \textit{et al.} \cite{Ahmed}  consider the minimization of waste and cutting patterns.   Table \ref{trabPCE} presents a summary of the literature  about the multi-objective CSP, highlighting the optimization criteria and the solution methods applied. Only recently the scalarization methods are being used to solve this problem. Additionally, the two-dimensional case has received considerably less attention.
 {\footnotesize
 \begin{table}[ht]
       \begin{center}
        \caption{Literature addressing the multi-objective CSP.}
 \begin{tabular}{ p{3cm}|p{5cm}| p{4cm} }
   \hline Author and Year  & Problem and objectives addressed & Solution Method\\
  \hline  W{\"a}scher (1990) \cite{wascher1990lp} & {1D}-CSP minimizing the total {number of objects, waste, overruns, leftovers}    & Five Phases Interactive procedure with local information\\
  \hline Sinuany-Stern and Weiner (1994) \cite{Sinuany-Stern-Weine} & {1D}-CSP minimizing  waste as a main optimization criterion and  leftovers as secondary criterion & Three Phases Heuristic Method \\
  \hline
   Kolen and Spieksma (2000) \cite{Kolen} & {1D}-CSP minimizing total waste  and  setup. &  \textit{Branch and Bound} Method\\
  \hline
 Respicio and Captivo (2005) \cite{Respicio}   &  {1D}-CSP
with the sequencing of the cutting patterns, minimizing the {maximum number of open stacks} and  \break {order spread}. & Evolutionary Algorithm\\
  \hline
Mu{\~ n}oz \textit{et al.} (2007)\cite{Munoz}   & {1D}-CSP minimizing of the {maximum number of open stacks}, total waste, and setup. & Genetic Algorithm \\
\hline  Gomez and Terashima-Maríns (2010)
\cite{Gomez}  &  {2D}-CSP (irregular)
 minimizing  {total waste } and  {total time to generate cutting patterns}. & Evolutionary Algorithm\\
 \hline  Salles-Neto \textit{et al.} (2014)
 \cite{Salles}   & {1D}-CSP 
 minimizing  leftovers,  setup and the total number of saw cycles. & Genetic Algorithm\\
  \hline Araujo \textit{et al.} (2014) 
  \cite{araujoPoldi2014}& {1D}-CSP minimizing total waste and  setup . & Genetic Algorithm\\
   \hline Arana-Jiménez and Salles-Neto (2017) \cite{Leduino} & {1D}-CSP minimizing the total number of {objects} and  setup. & Sufficient conditions for partial efficient solutions that might be used in heuristic methods.\\
  \hline Toscano \textit{et al.} (2017)\cite{Toscano}  &  {2D}-CSP minimizing the total number of objects and {saw cycles}. & Cutting plane Heuristic\\
 \hline Rangel and Sáez-Aguado (2017)\cite{Rangel2}  &  {1D}-CSP minimizing the total number of {objects} and {saw cycles}. & Scalarization methods\\
\hline Aliano-Filho \textit{et al.} (2018)  \cite{Angelo2018}& {1D}-CSP minimizing total waste and the total number of {cutting patterns}. & Scalarization methods\\
  \hline Mellouli \textit{et al.} (2019)
 \cite{Ahmed}  &  {2D}-CSP minimizing total waste and the total number of {cutting patterns}. & Genetic Algorithm\\
 \hline  Aliano-Filho \textit{et al.} (2021) \cite{Angelo2021}  & {1D}-CSP minimizing total waste and the total number of {cutting patterns}. & Scalarization methods\\
  \hline { Pierini and Poldi (2022) \cite{PoldiMO} }&  {1D}-CSP ( multiperiod) minimizing total waste and total inventory cost of objects and items. & Scalarization methods\\
  \hline
  \end{tabular}\\
 \label{trabPCE}
 \end{center}
  \end{table}
}

The focus of this article is to present a study of the bi-objective CSP (BOCSP) taking into account the minimization
of objects and saw cycles as well as the minimization
of objects and setup. Both the 1D case and the 2D regular case are studied. The main contributions of this research are: the
mathematical optimization model used to represent the problem;  the use of
dynamic column generation embedded into scalarization methods to solve the BOCSP; a computational study that highlights the trade-off among the total number of objects and saw cycles according to the variation of different parameters of the problem; reduction of the gap in the literature related to the BOCSP, in particular the two-dimensional case. Most of the scalarization methods employed to solve the multi-objective CSP up to now consider that the set of cutting patterns are generated 
\textit{a priori}. The methods proposed here consider that they are dynamically generated. To the best of our knowledge, besides us, only \cite{PoldiMO} has employed  such method  to solve the BOCSP. 

The remainder of this text is organized as follows. In Section
\ref{mod} the  problem and the mathematical model used to represent it are described and in Section \ref{boco} some basic concepts of multi-objective optimization  as well a brief description of the  three scalarization methods employed
are reviewed. In Section \ref{BOCG} the algorithms proposed to solve the BOCSP for
the one-dimensional case  (1D-BOCSP) and for the bi-dimensional case (2D-BOCSP) are detailed. The computational study design as well as the analysis of its results are presented in Section \ref{comp}. Final considerations are given in Section \ref{concl}.
  
\section{Problem Description and Model formulation}\label{mod}
The cutting stock problem considered here can be stated
as defining how to cut a  set of single type large pieces (denoted by
objects) of thickness $t$ to produce
a set of $m$ smaller pieces (denoted items) with the same thickness $t$ in order to
fulfill a demand $d_i$ for each item $i$, $i=1...m$ according to a given optimization criteria. We consider that the objects are available in stock in enough quantity and that the cutting machine can be used  to cut many objects, all at once, by just stacking them one on top of the other. The maximum number $p$ of objects that can be cut at a time (saw capacity) is determined by the thickness of each object  and by the saw height,
$h$, and is computed according to (\ref{maxcut}).  A saw cycle stands for all the cutting machine operations to cut one object (or a stack of objects simultaneously). When the number of objects in a saw cycle equals $p$, the cycle is named \textit{complete saw cycle}.
\begin{equation}
    p=\left\lfloor \frac{h}{t} \right \rfloor. \label{maxcut}
\end{equation}

Some elements must be taken into account when defining the optimization criteria of the CSP addressed here. Lets consider the mono-objective CSP with the usual criterion of minimizing waste computed as minimizing the total number of objects (CSP-O). As pointed out in \cite{Toscano}, if the diversity of the cutting patterns in the optimal solution is reduced, the frequency of a given cutting pattern might increase (on average) which might contribute for reducing the total number of saw cycles. So, using the setup minimization as an optimization criterion might contribute to reduce the total number of saw cycles, although this is not guaranteed. Another interesting feature might occur if the frequency of each one of the cutting patterns of the optimal solution is a multiple of $p$. In that case, as showed in \cite{Yanasse2008}, the mono-objective CSP considering the minimization of the total number of saw cycles is equivalent to the mono-objective CSP minimizing the total number of objects with the demands scaled by $p$, except for a constant. However, this is unlikely to happen, unless $p = 1$. 

At scenarios of high demand, when the cost of the machine time might dominate the other costs of the production, a desired decision  criterion  is to reduce the total number of saw cycles, assuming that the machine takes about the same time to cut one or $p$ objects. Still, it can conflict with the criterion  of minimizing the total number of objects, since in case of excess of demand being accepted, a better use of the saw capacity (complete saw cycles) might imply in an increase of the number of objects cut. So there is interest in finding a compromise between these two objectives and a multicriteria point of view should be explicitly employed. This aspect is particularly important in the furniture industry \cite{Toscano}.

From now on the bi-objective CSP problem addressed here  considering the minimization of the total number of objects and of saw cycles is referred to as BOCSP. An interesting feature of this problem occurs when the saw capacity $p$ is greater or equal to $d_{max} = max\{d_i, i=1...m\}$.  In that case, the   BOCSP is equivalent to the bi-objective CSP problem considering the minimization of the total number of objects and of Setups (BPCSP-S) \cite{Salles}. This feature is explored in the computational study presented in Section \ref{comp}.

To build a mathematical optimization model to represent the BOCSP, suppose that there are $n$ possible cutting patterns that have been
generated \emph{a priori}. The cutting patterns can be represented by $A_j=[a_{ij}]$,
a $m$-dimensional column vector with each element $a_{ij}$ being the
number of items $i$ in pattern $j$, $j=1,...,n$. Let $x_{j}$ and
$y_j$ be the number of objects and cycles associated with the cutting
pattern $A_j$, respectively. Letting:
\begin{eqnarray}
\label{MF1}f_1(x,y)&=&\sum_{j=1}^{n}x_{j} \mbox{ be the total number of objects}, \\
\label{MF2}f_2(x,y)&=&\sum_{j=1}^{n}y_{j} \mbox{ be the total number of  cycles},
\end{eqnarray}
the bi-objective mathematical optimization model (BMO) is given by
(\ref{M1})-(\ref{M6}) and is a variation of the model presented
in \cite{Yanasse2008}.
\begin{eqnarray}
{min} &   \quad (f_1(x,y), f_2(x,y)).& \label{M1}\\
\nonumber \mbox{subject to}& & \\
&{\sum_{j=1}^{n}a_{ij}x_{j}}\geq d_i, & i=1,...,m. \label{M2}\\
&{\sum_{j=1}^{n}a_{ij}y_{j}}\geq \lceil{d_i/ p}\rceil,
&i=1...m \label{M4}\\
&{x_j \leq py_j},         &j=1,\ldots,n.\label{M3} \\
&  x_{j}, y_j  \in \mathbb{Z}^+ & j=1,\ldots,n.\label{M6}
\end{eqnarray}
The constraints (\ref{M2}) guarantee that the demand is satisfied
with overproduction and the constraints (\ref{M3}) control the
number of cycles associated to each cutting pattern. The computational study presented in
\cite{Rangel} shows that including the redundant  constraints (\ref{M4})
in the model improves the solution process of the mono-objective
version of the model that minimizes ($f_1(x,y)+f_2(x,y)$).

The model (BMO)  has $2n$ columns and ($2m+n$) constraints. The value of $n$ can grow exponentially according to the value of $m$ and the
items dimensions. A common approach used to cope with the high number of columns (\textit{\textit{e.g.}} \cite{Angelo2021}) is to use a restricted version of the model in which  only a subset of cutting patterns (columns) is considered, this approach is named here as Static Column Generation (SCG). The other way is to generate the columns as they are needed, named here as Dynamic Column Generation (DCG) and is further discussed in Section \ref{secGCd}. Either way, a initial set of columns should be generated before a scalarization method is employed to solve instances of  the (BOCSP).  

The model (\ref{M1})-(\ref{M6})  can the used to represent both the one-dimensional (1D-BOCSP) and the two-dimensional (2D-BOCSP) cases of the (BOCSP),the difference lies in the procedure used to generate the columns (cutting patterns) for each case. This is further discussed in the section \ref{SecGP}.

\subsection{Procedures to generate 1D  and  2D cutting patterns }\label{SecGP}

Before going into more details about the cutting pattern generation for the BOCSP  it is necessary to define  the number of dimensions that are relevant to the cutting procedure as well as what are constraints related to the cutting equipment. In the context of this research, for the 1D-BOCSP only one dimension is relevant for the cutting procedure, both  for the object and for the items, which is their length, $L$ for the objects and $l_i$, $i=1...n$, for the items. There is no constraint related to the cutting procedure.

\begin{figure}[ht]
\center
\subfigure[PD1][A two-dimensional cutting pattern]{\includegraphics[width=3.0cm]{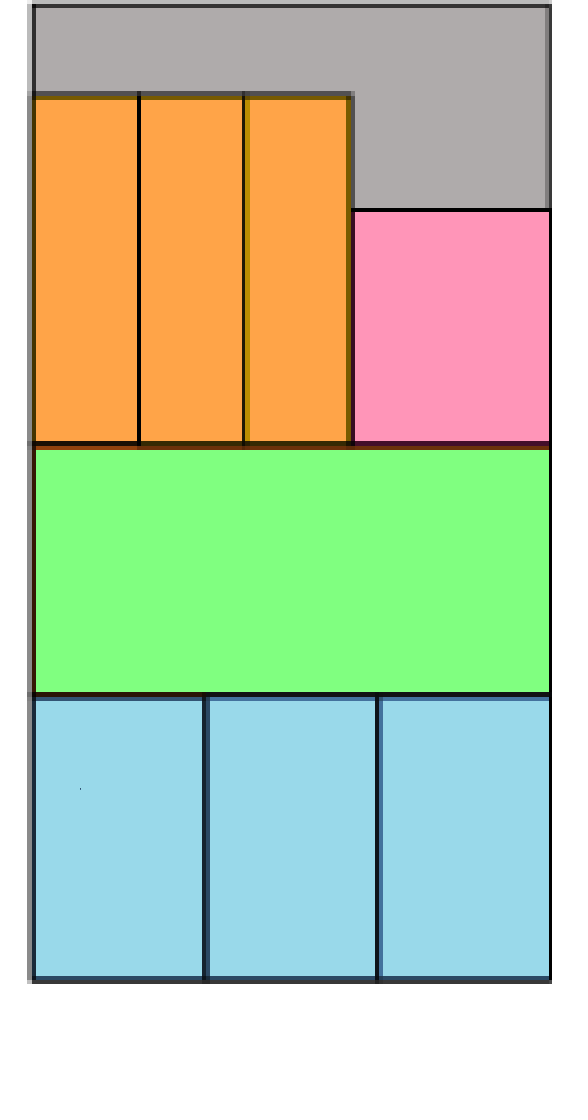}\label{pd1}}
\qquad
\subfigure[Est1][1st stage]{\includegraphics[width=3.0cm]{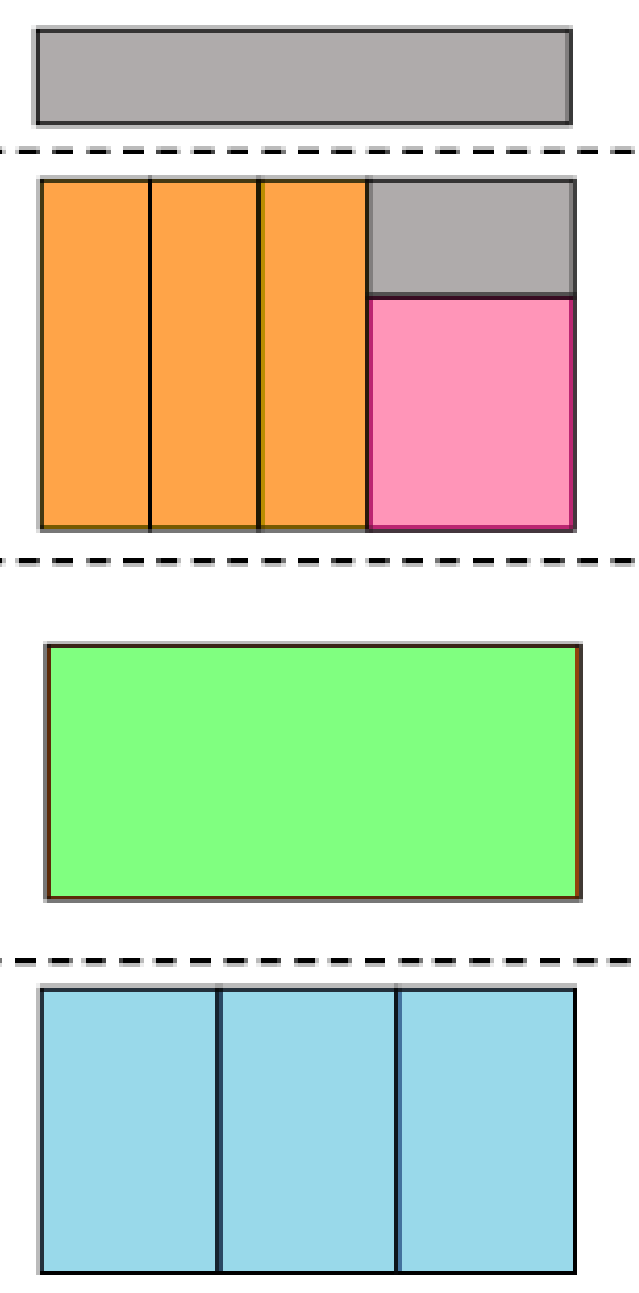}\label{ESt1}}
\qquad
\subfigure[Est1][2nd stage]{\includegraphics[width=3.0cm]{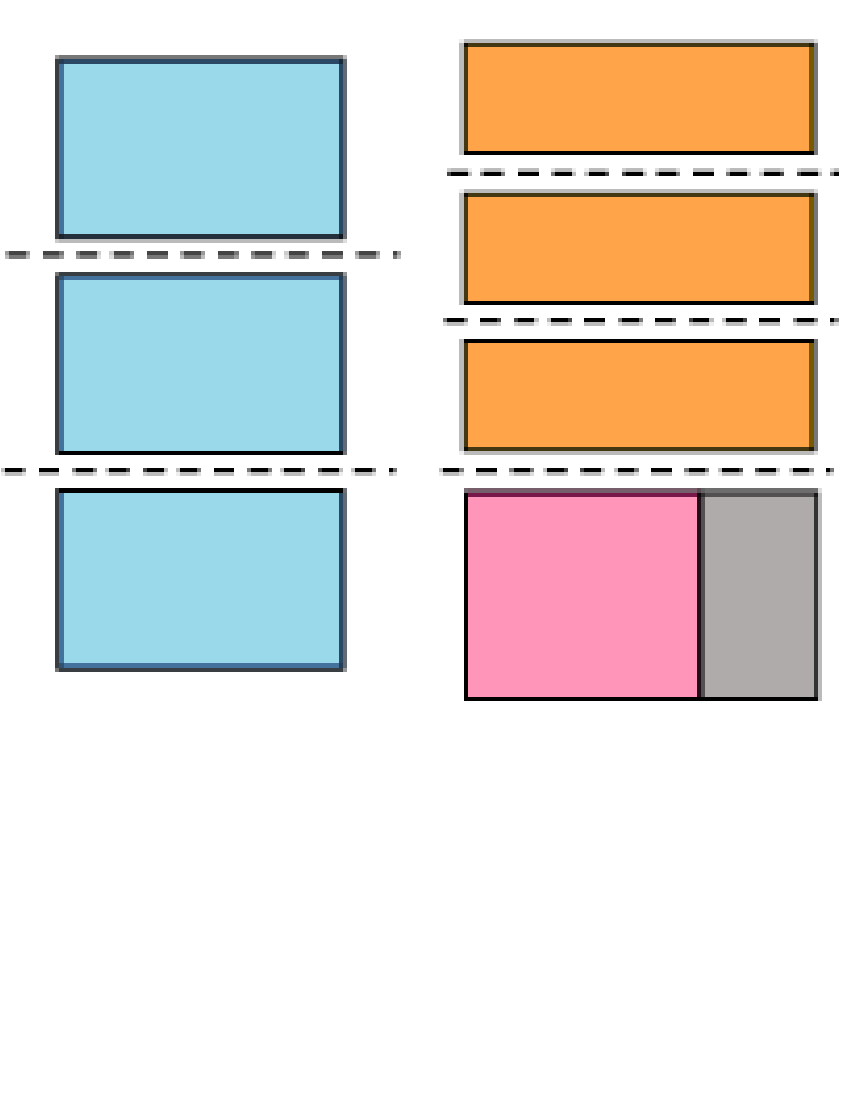}\label{ESt2}}
\caption{A two-stage orthogonal guillotine non-exact cutting pattern.}
\label{CPGO}
\end{figure}

As for the 2D-BOCSP, two dimensions are relevant for the cutting procedure, the length $L$ and width $W$ for the object,  and the lengh $l_i$, $i=1...n$, and width $w_i$, $i=1...n$, for the items. Due to constraints related to the cutting equipment, it is considered that only orthogonal guillotine cuts, that is, cuts made from one side of the object to the other (guillotine) and parallel to the other side of the object (orthogonal) are allowed. Also, the object is rotated only once to obtain the items (two-stage cut).  Figure \ref{pd1} shows an example of a two-stage orthogonal guillotine cutting pattern. The cutting procedure goes as follows. First the object is divided into strips (Figure \ref{ESt1}). Those strips are then rotated by 90 degrees to be further cut in order to obtain the items (Figure \ref{ESt2}). Since not all items are obtained with a single rotation, this cutting pattern is classified as  non-exact. To obtain some of the items, an extra rotation is necessary, or they are sent to a secondary cutting machine.

To generate a cutting pattern we have to solve a Single Large Object Placement Problem (SLOPP) \cite{Wascher2007}. A feasible solution to SLOPP represents the way the items are positioned (allocated) in the object before cutting. As with other Cutting and Packing Problems (C\&P), two basic conditions must be satisfied to obtain a feasible solution: the  items included cannot exceed the dimensions of the object (containment constraint) and they cannot overlap (non-overlapping constraint). In \cite{scheithauer2017} a general model for (C\&P) is proposed considering that a value is assigned to each item and that the optimization criterion is to maximise the total value of the elements included in the object. 

In the context of the column generation procedure, the problem solved to generate a column is named \textit{pricing subproblem} since it is related to the pricing operation in the simplex method (\textit{e.g.} \cite{LivroCG}). The items values are the dual variables associated with the demand constraints (given by either (\ref{M2}) or (\ref{M4})) in the optimal solution of the linear relaxation of (\ref{M1})-(\ref{M6}). To formulate the subpricing's model for the 1D case, let's assume that the items will be placed in the object one after the other, without overlapping, and such that the sum of the lengths of the selected items does not exceed the length of the object ($L$). In this case, the containment and non-overlapping constraints can be expressed by a single mathematical expression, the inequality (\ref{pmoch2}), in which $a_{i} \in Z$ represents the number of items of length $l_i$ included in the object (cutting pattern), $i=1...m$. Assuming that the objective is to determine a subset of items that provides the highest value, and that $\pi_i$ is the value of item $i, i=1...m$, the 1D \textit{pricing subproblem}  can be  represented by the model  \eqref{pmoch1}-\eqref{pmoch3}. This model is also used to represent another packing problem known as the Integer Knapsack problem (\textit{\textit{e.g.}} \cite{Martello1990}). 
\begin{eqnarray}
&max&  \sum ^{m}_{i=1}\pi_i a_{i} \label{pmoch1} \\
&s.t.  &\sum ^{m}_{i=1}l_j \cdot a_{i}\leq L, \label{pmoch2} \\
& &a_{i} \in \mathbb{Z}_{+}. \label{pmoch3} 
\end{eqnarray}
 
Several approaches have been proposed to solve the 2D  version of the (SLOPP) including  mathematical optimization models solved by general purpose solvers and/or matheuristics (\textit{e.g.} \cite{Lodi, Malaguti, Mateus2020}). Computational studies have shown that a good strategy to generate a two-stage orthogonal guillotine non-exact cutting pattern is to use the mathematical optimization model proposed in \cite{Lodi}. In the model named (M1), the items are all considered distinct, so that the parameter $m$ is  used to define the number of types of items according to the number of different dimensions $l_i\times w_i$, and the parameter $s$ for the total quantity of items demanded, that is, $s= \sum ^{m}_{i=1}d_i$, where $d_i$ is the demand for the item $i$. The items are sorted in non-ascending order according to their width $w_i$ ($w_1\geq w_2\geq \ldots \geq w_m$). In the model, it is assumed that $s$ strips can be initialized, one  for each width $w_k$, $k\in \{1, \ldots , s\}$. A binary decision variable $x_{ik}$ is defined to represent if the item $i$ is allocated in strip $k$ ($x_{ik}=1$), or not ($x_{ik}=0$). The model M1 is given by \eqref{LM1}-\eqref{LM5}.
\begin{eqnarray}
&max&  \sum ^{s}_{i=1}\pi_i\sum ^{i}_{k=1} x_{ik}\label{LM1} \\
& s.t.  &\sum ^{i}_{k=1}x_{ik}\leq 1, \quad \quad \quad \quad \quad \quad i=1, \ldots, s \label{LM2} \\
& &  \sum ^{s}_{i=k+1}l_i x_{ik} \leq (L-l_k)x_{kk}, \quad  \quad k=1, \ldots, s-1 \label{LM3} \\
& &   \sum ^{s}_{k=1}w_k x_{kk} \leq W,\label{LM4} \\
& &x_{ik}\in \{0,1\}, \quad \quad \quad \quad \quad i=1, \ldots , s \quad k=1, \ldots , s;\label{LM5} 
\end{eqnarray}
in which $\pi_i$ is the value of the item $i$, and the objective function \eqref{LM1} maximizes the total  value of the items allocated in the object. The constraints \eqref{LM2} guarantee that each item is allocated only once and in the strip with width greater or equal to its width. The constraints \eqref{LM3}, similar to \eqref{pmoch2}, guarantees that the total length of the items allocated in strip $k$ does not exceed the strip length $L$. The constraint  \eqref{LM4} guarantee that total width of the strips allocated in the object do not exceed the width of the object, and the constraints \eqref{LM5} indicate the domain of the variables. 

To consider the case that rotation of the items  is allowed, for each item $i$ a complementary item $i+s$ is created so that, $l_{i+s}=w_i$ and $w_{i+ s} = l_i$.  Model M1 with rotation of items  (M1-rot)  is obtained replacing the constraints \eqref{LM2} with \eqref{LMrot} and $s$ with $2s$, considering that all the items can be rotated.
\begin{equation}\label{LMrot}
    \sum ^{i}_{k=1}x_{ik}+\sum ^{\theta _i}_{k=1}x_{\theta_{i}k}\leq 1 \quad i=1, \ldots , 2s \mbox{ and } i<\theta _i,
\end{equation}
where $\theta _i$ indicates the rotated item corresponding to item $i$. 
The number of variables in the model is directly related to the
number of copies created from the demand for each item. Therefore to reduce the total number of variables, the demand fulfillment is controlled through Equation \eqref{dbar} and taking $s=\sum ^{m}_{i=1}\bar{d_i}$.
\begin{equation}\label{dbar}
\bar{d_i}= min\left\{d_i,\left\lfloor\frac{L\cdot W}{l_i\cdot w_i}\right\rfloor \right\}, i=1, \ldots , m.
\end{equation}
A trivial set of cutting patterns can be generated considering that a single item is allocated at its maximum and is classified as  ``maximal homogeneous" or simply ``homogeneous" . For the one-dimensional case the set with $m$ homogeneous cutting patterns can be computed  according to \eqref{HCP-1D}  and for the two-dimensional case with rotation the set with up to $2m$ homogeneous cutting patterns can be computed according to  \eqref{HCP-2D}.
\begin{equation}\label{HCP-1D}
a_{ij}= \lfloor\frac{L}{l_i}\rfloor,  \mbox{ if } i=j; \quad 0, \mbox{otherwise} .
\end{equation}
\begin{equation}\label{HCP-2D}
a_{ij}= \lfloor\frac{L}{l_i}\rfloor \cdot \lfloor\frac{W}{w_i}\rfloor,  \mbox{ if } i=j;   \quad 0, \mbox{otherwise}.
\end{equation}

\section{Basic concepts and  approaches to solve the BOCSP}\label{boco}
The model (BMO) can be restated according to (\ref{biobj})-(\ref{bir}).
\begin{eqnarray}
min &     & (f_1(x,y), f_2(x,y)) \label{biobj} \\
s.t.&      & (x,y)\in \mathcal{X} \label{bir},
\end{eqnarray}
with the decision space given by $\mathcal{X}= \{(x,y):$ satisfy simultaneously $ (\ref{M2})-(\ref{M6})\}$ and the criteria space defined as
$\mathcal{Z} = \{z=(f_1(x,y), f_2(x,y)) : (x,y) \in \mathcal{X} \}$.

Since, generally, no solution $(x,y)\in \mathcal{X}$ simultaneously minimizes both
objectives, one searches for solutions that accomplish an acceptable trade-off between the objectives. The bi-objective problem is considered to be solved when the set of efficient or Pareto optimal points has been found \cite{Ehrgott2003}. For ease of explanation, first the main definitions that are used in the remainder of the text are presented. 
\begin{itemize}
\item $(x^1,y^1) \in \mathcal{X}$ dominates $(x^2,y^2) \in \mathcal{X}$ if $f_i(x^1,y^1) < f_i (x^2,y^2)$, $i=1,2$, and if an equality holds for only one $i \in \{1,2\}$ it is defined as weak dominance;
  \item $(x^*,y^*)$ is said to be an efficient solution (or Pareto optimal solution, or nondominated solution) if does not exist a solution $(x,y)$ $\in$ $\mathcal{X}$ that dominates $(x^*,y^*)$ (In case of the weak dominance it is said weakly Pareto optimal solution);
  \item $\mathcal{X}^*\subseteq \mathcal{X}$ denotes the set of all efficient solutions in $\mathcal{X}$;
   \item $z^*=(f_1(x^*,y^*), f_2(x^*, y^*))$ is said to be a nondominated point (nondominated vector) in the criteria space if $(x^*,y^*)$ $\in$ $\mathcal{X}^*$;
  \item $\mathcal{Z}^*$ denotes the set of all nondominated points, also called Pareto front, and $N=|\mathcal{Z}^*|$ denotes the cardinality of the set $\mathcal{Z}^*$;
  \item A nondominated point $\bar{z} \in  \mathcal{Z^*}$ is called non-supported if it is dominated by a convex combination (which does not belong to  $\mathcal {Z^*}$) of
other nondominated  points. Otherwise it is said to be
supported;
\item $z^{I}=(f_1^{I}$, $f_2^{I})$ is said to be \textbf{Ideal Vector} if $f_1^I = min\{f_1(x,y)\mid (x,y) \in \mathcal{X}\}$ and  $f_2^I = min\{f_2(x,y)\mid (x,y) \in \mathcal{X}\}$;
\item Let $z^{1}=(f_1^1,f_2^1) \in Z$ and $z^{2}=(f_1^2,f_2^2) \in Z$ with $f_1^1 \leq f_1^2$ and $f_2^2 \leq f_2^1)$. $R(z^{1},z^{2})$ denotes the rectangle defined by the points $z^{1}$ and $z^{2}$;
\item $z^{L1}$ and $z^{L2}$ are said to be \textbf{Lexicographic Vectors} if $z^{L1}=(f_1^{L1}$, $f_2^{L1})$ and $z^{L2}=(f_1^{L2}$, $f_2^{L2})$ where:\\ $f_1^{L1}=f_1(x^{*1},y^{*1}) = min\{f_1(x,y)\mid (x,y) \in \mathcal{X}\}$; \\ 
$f_2^{L1} =f_2(x^{*1},y^{*1})=min\{f_2(x,y)\mid (x,y) \in \mathcal{X}$ and $ f_1(x,y) = f_1^{L1} \}$; \\ 
$f_2^{L2}=f_2(x^{*2},y^{*2}) = min\{f_2(x,y)\mid (x,y) \in \mathcal{X}\}$; \\ 
$f_1^{L2} =f_1(x^{*2},y^{*2})=min\{f_1(x,y)\mid (x,y) \in \mathcal{X}$ and $ f_2(x,y) = f_2^{L2} \}$.
\end{itemize} 

In the context of multi-objective optimization, the solution methods are classified into three categories: \textit{a priori}, \textit{a posteriori} and \textit{interactive}, according to the role of the decision makers in the problem solving process. The \textit{a priori} methods are those in which the decision makers provide information and preferences regarding the problem prior to the resolution process. In \textit{a posteriori} methods, the set of all efficient solutions is generated and the decision makers analyze the set according to their preferences. \textit{Interactive} methods are those in which the decision makers introduce preferences in an interactive way during the solution process. This research focuses on \textit{a posteriori} methods aiming to obtain the Pareto front or  a good approximation of it. In this solution strategy, scalarization methods are used, which transform a multi-objective problem into a sequence of mono-objective problems, which under certain conditions, generate efficient solutions.

In Section \ref{BOCG} a discussion of column generation embedded into the scalarization method is presented followed by adaptations of the three methods chosen to solve the BOCSP. The  Lexicographic $\epsilon$-Constraint method ($LEC$) (\textit{\textit{e.g.}} \cite{Saez}) was chosen for being a classic variation of the $\epsilon$-Constraint method, and in order to assess whether its performance is impaired by the need to solve two subproblems at each iteration. The  Frontier Partitioner Algorithm $FPA$ \cite{Santis2020} was chosen for its ease of implementation and the fact that it has a better performance when compared to the \textit{Balanced Box} method \cite{Boland2015}. Finally, the  Augmented Weighted Tchebycheff Method ($AWT$)  was chosen based on the discussions presented in \cite{Angelo2018, Angelo2021} that indicate its reduced computational time to obtain an approximation of the Pareto front.
To the best of our knowledge, there are no studies applied to (BOCSP) that embeds dynamic column generation into scalarization methods and that compare these three methods. 

\section{ Scalarization methods for the BOCSP (1D and 2D)}\label{BOCG} 
One of the main contributions of this research is the study of the integration of the column generation technique \cite{LivroCG} to the scalarization methods  $LEC$,  $FPA$ and $AWT$ to solve the BOCSP (1D and 2D). Basically, these methods consist of transforming the solution of a multi-objective problem into the solution of a sequence of mono-objective problems and thus finding the Pareto front (or an approximation to it). Each of these mono-objective problems is called \textit{Internal Scalarization Problem} (ISP). In Section  \ref{secGCd} we present the multi-objective column generation method and in Sections \ref{LER}-\ref{FPA} the scalarizations applied to solve the BOCSP.

\subsection{Column Generation for the BOCSP}\label{secGCd}
As discussed in Section \ref{sec:intro}, the Column Generation (CG) technique has not been much employed to solve multi-objective optimization problems. The study presented in \cite{Artigues2018} shows that the few CG techniques that are found in the literature to solve multi-objective problems follow a basic structure that can be summarized as follows. First, transform the multi-objective problem into a mono-objective problem using some scalarization method and considering a subset of columns generated \textit{a priori}, this problem is called Restricted Master Problem (RMP). Next, the CG is applied to the linear relaxation of the subproblem generated for one or more  parameters of the scalarization (LRMP). If new columns are added to the LRMP the process is repeated, otherwise it is finished. If the CG procedure is applied to solve the LRMP for all the parameters of the scalarization method it is called  ``point-by-point search'' \cite{Artigues2018}.

The algorithm implemented in this study is called Dynamic Column Generation (DCG) and consists in three main steps:
\begin{itemize}
\item Generate an initial  subset of columns to compose the RMP ($CG_a$).
\item Compute the lexicographic points (See Section \ref{boco}) applying the CG algorithm (Algorithm \ref{AlgGC}) to solve each one of the  associated optimization problems.
\item Apply the ``point-by-point search'' to solve the  BOCSP.
\end{itemize} 
The inicital subset of columns are generated \textit{a priori } ($CG_a$) 
applying twice the Algorithm \ref{AlgGC}  \cite{Gilmore1961, Rangel2} to mono-objective versions of the problem, one for each objective function.  The first  subset of columns is obtained by solving the CSP considering the minimization
of the total number of objects, given by (\ref{MF1}), and the constraints defined by (\ref{M2}). Then the second subset of columns is
obtained by solving the CSP considering the minimization of the total number of cycles, given by (\ref{MF2}), and the constraints
defined by (\ref{M4}). The associated matrices are them merged eliminating the duplicated columns. The initial set of columns for both problems is composed by homogeneous cutting patterns. The LRMP problem is obtained by redefining  the domain of the variables $x,z \in  \mathbb{Z}_{+}^n$  to $x,z\in \mathbb{R}_{+}^n$.

In  the CG procedure it is necessary to solve a pricing subproblem (\textit{subpricing}) to obtain a new column (new cutting pattern). For the 1D case the  \textit{subpricing} is the Knapsack problem and for the 2D case the  \textit{subpricing} is the problem M1-rot, both presented in Section  \ref{SecGP}. 

\begin{center}
\begin{minipage}{15cm}
{\footnotesize
\begin{algorithm}[H] \label{GC}
   \SetKwData{In}{\textbf{Input}} 
   \SetKwData{Out}{\textbf{Output}}
   \SetKwData{Calc}{\textbf{Compute}}
   \SetKwData{Resol}{\textbf{Solve}}
   \SetKwData{Atual}{\textbf{Update}}
   \SetKwData{Const}{\textbf{Build}}
   \SetKwData{Ret}{\textbf{Go back}}
   \SetKwData{Faca}{\textbf{Do}}
   \SetKwData{Par}{\textbf{Stop!}}
   \SetKwData{Cha}{\textbf{Call}}
    \SetKwData{Ins}{\textbf{Insert}}
   \SetAlgoLined
  \In BOCSP data\text{,} $A$ (matrix with an initial set of cutting patterns)\;
  \Out $A$ matrix updated and the optimal solution to the associated LRMP\;
	 \Resol LRMP associated with the current matrix $A$ and get the dual variables associated with the demand constraints\;\label{reopt}
	 \Resol \textit{subpricing} to get a new column (cutting pattern) and the associated reduced cost\;
	\eIf{$\mbox{ Reduced cost}\geq 0$}{ \Par The current solution of LRMP is optimal}{\Ins the new cutting pattern in matrix $A$ and reoptimize (\Ret line \ref{reopt})}	             
	\caption{Column Generation (CG)}\label{AlgGC}
\end{algorithm}
}
\end{minipage}
\end{center}

\subsection{The lexicographic  $\epsilon$-constraint method} \label{LER}
The $\epsilon$-constraint  method  consists of  obtaining the set $\mathcal{Z}^*$ solving the subproblem (\ref{Sub1}) for different values of the parameter $\epsilon$. In the lexicographic $\epsilon$-constraint method ($LEC$)  (\emph{\textit{e.g.}} \cite{Saez}) the set $\mathcal{Z}^*$ is obtained by iteratively solving two subproblems. For the BOCSP, the first subproblem ($Sub1$ as stated in (\ref{Sub1})) aims at minimizing the total number of objects ($f_1(x,y)$) and the second subproblem ($Sub2$ as stated in (\ref{Sub2})) aims at minimizing the total number of cycles ($f_2(x,y)$). The $\epsilon$-constraint (\ref{c_eps1}) is included in the definition of subproblem $Sub1$ and an adjustment constraint (\ref{c_adj}) in the definition of the subproblem $Sub2$. At each iteration $t$ the right
hand side  of (\ref{c_eps1}) is updated according to the value of $\epsilon=(\widehat{f}_2^{t-1}-1)$ obtained with the solution of $Sub2$ in the iteration ($t-1$) and the right-hand side of (\ref{c_adj}) is updated according to the value of ${f}_1=\widehat{f}_1^{t}$ obtained with the solution of $Sub1$. The algorithm halts when the current $Sub1$ is infeasible.

\begin{eqnarray}
\label{Sub1}\mbox{Sub}1:\quad \widehat{f}_1= \{ min \quad f_1(x,y):
(x,y) \in \mathcal{X} \cap
\{(\ref{c_eps1})\}\}\\
\label{Sub2}\mbox{Sub}2:\quad \widehat{f}_2= \{ min \quad f_2(x,y):
(x,y) \in \mathcal{X} \cap \{(\ref{c_adj})\}\}.
\end{eqnarray}

\begin{eqnarray}
\label{c_eps1}\sum_{j=1}^{n}y_{j}\leq \epsilon&\\
\label{c_adj}\sum_{j=1}^{n}x_{j}\leq f_1
\end{eqnarray}

The solution of Sub2 avoids the necessity of eliminating the
nondominated solutions and the nondominated points from the sets
$\mathcal{X}^*$ and $\mathcal{Z}^*$. Among all the feasible solutions with the
number of objects being equal or less than equal to the objective
value of Sub1, one searches for the ones with the minimal possible
number of cycles. If Dynamic Column Generation (DCG) is used, the Algorithm \ref{AlgGC} is applied to solve $Sub1$ and $Sub2$ for each value of the parameter $\epsilon$. Otherwise, all the subproblems are solved considering the same set of columns (set $CG_a$ as described in Section \ref{secGCd}), the Static Column Generation method (SCG).

\subsection{The  Frontier Partitioner Algorithm}\label{descFPA}\label{FPA}

The Frontier Partitioner Algorithm ($FPA$) proposed in \cite{Santis2020} is a \textit{Branch and Cut} algorithm that consists of introducing inequalities  inducing a partition of the criterion space and using a scalarization method to obtain a nondominated point. Considering the existence of the ideal vector $z^{I}_i, i=1, 2$ and also that for each of the objective functions, there is a positive value that is greater than the distance between the image of two feasible solutions (\textit{i.e.} the objective functions are $\gamma$-positive) it is possible to show that the Pareto front associated with the problem is finite. Furthermore, if the $ISP$ have an optimal solution or it is possible to show that it is infeasible,  the method is capable of producing the entire Pareto front $\mathcal{Z}^{*}$. For the BOCSP studied here, $f_1$ and $f_2$ are $\gamma$-positive, for $\gamma \geq 1$.

The scalarization method used is a combination of the Weighted Sum method and the Lexicographic method called \textit{Custom Weighted-sum Scalarization} (CWS). Considering a permutation $(i_{1},i_{2})$ of the set $\{1,2\}$, let $(\hat{x},\hat{y})$ be a solution of the associated lexicographic problem and $(\bar{x},\bar{y})$ a solution of the lexicographic problem associated with the reverse permutation. The respective points $\hat{z}$ and $\bar{z}$ are nondominated and different (otherwise the ideal vector would be optimal to the BOCSP). If $\gamma={min}_{i=1,2}\{\gamma_i\}>0$ and $\zeta \in (0,\gamma)$, the weights $w(\zeta)\in \mathbb{R}^2_{>}$ calculated according to \eqref{defw} are valid for the CWS and allow the scalarization \eqref{SPPer} to be applied in the branching scheme as the weights are valid for the descendant nodes.
 \begin{equation}\label{defw}
w(\zeta)_i=
\left\{\begin{array}{cc}
\dfrac{\gamma -\zeta}{\hat{z}_{i_1}-\bar{z}_{i_1}}, & \mbox{ se } i=i_1 \\
1, & \mbox{ se } i=i_2
\end{array}\right.
\end{equation}
\begin{eqnarray}
& min & \quad F(x,y)=w(\zeta)_1f_{1}(x,y)+w(\zeta)_2f_{2}(x,y)\label{SPPer} \\
     & s.t. & \quad (x,y)\in \mathcal{X}.\nonumber  
\end{eqnarray}
\vspace{-0.70cm}

The FPA  starts from a nondominated point and defines two subproblems so that the chosen nondominated point and all points that are dominated by it are infeasible for the two subproblems. Thus, consider a generic \textit{node} $k$ from the tree, the subproblem ($S^k$) is represented in \eqref{POBIk2}.
  \begin{eqnarray}
(S^k)\quad & min & \quad F(x,y)\label{POBIk2} \\
      & s.t. & \quad (x,y)\in \mathcal{X}^k,\nonumber
    \end{eqnarray}
in which $\mathcal{X}^k\subseteq \mathcal{X}$ is the set obtained from the intersection of $\mathcal{X}$ with one of the inequalities defined in \eqref{eqram}. For $k=0$, (\textit{root node}) we have ($S^0$)=BOCSP and $\mathcal{X}^0=\mathcal{X}$.

If the subproblem $S^k$ has an optimal solution, the feasible region associated to \textit{node} $k$ can be partitioned generating two \textit{children nodes}. Let $(\hat{x},\hat{y})^k \in \mathcal{X}^k$ be the optimal solution of $S^{k}$ and $\hat{z}^k=F( \hat{x},\hat{y})$. Taking $\epsilon_i\in (0,\gamma_i], i=1, 2$ and considering the inequalities \eqref{eqram}, the subproblems associated with \textit{child node} $i$, $i=1, 2$, is given by \eqref{POBIk1}.
\begin{eqnarray}
&&f_i(x,y)\leq \hat{z}_i -\epsilon_i , i=1, 2. \label{eqram}
\end{eqnarray}
\begin{eqnarray}
& min & \quad F(x,y) \nonumber \\
     & s.t. & \quad (x,y)\in \mathcal{X}^{k}_{i}, \label{POBIk1} 
\end{eqnarray}
 with   $\mathcal{X}^{k}_{i} = \mathcal{X}^{k} \cap \{x\in \mathbb{R}^{n}: f_i(x)\leq \hat{z}^{k}_{i} - \epsilon_{i}\}$, $i\in \{1,2\}$.
 
 To analyze the convergence of the FPA, the authors show that the inequalities (\ref{eqram}) used to generate the \textit{ children nodes} define a partition in the decision space and that each \textit{node} examined can be pruned by infeasibility or produces a nondominated point not yet found. Applying the $CWS$ scalarization, the number of subproblems is reduced to $|\mathcal{Z}^{*}|+1$ (best bound found so far) because it is not necessary to solve one of the subproblems generated in each \textit{node}. If the permutation $\{1,2\} $ ($\{2,1\}$) is used the subproblem associated with $ \mathcal{X}^{k}_2$ ($ \mathcal{X}^{k}_1$) is infeasible. 

 In the case of the $BOCSP$, we can take $\gamma_i=1, i=1,2$ in the calculation of $w(\zeta)_i$ since the two objective functions are $\gamma$-positive for $\gamma \geq 1$. From the value of $\gamma$, it is necessary to choose a value for $\zeta_i\in (0,\gamma_i)$ and a value for $\epsilon_i\in (o,\gamma_i]$. Thus, knowing $(\bar{x},\bar{y})$ the ISP solutions in iteration $k$ such that $\bar{z}^k=(f_1(\bar{x}),f_2(\bar{y}))$, the subproblem to be solved in the iteration ($k+1$) is given by \eqref{POBIk}- \eqref{mFPA4a}. The $FPA$ algorithm, including the possibility of column  generation when solving each \textit{ISP}, is displayed in Algorithm \ref{algFPA}.
\begin{eqnarray}
(S^{k+1})\quad & min & \quad F(x,y)\label{POBIk} \\
      & s.t. & \quad (x,y)\in \mathcal{X}^k \cap \{\eqref{mFPA4a}\}. \label{mFPA4}
\end{eqnarray}
\begin{equation}
  f_{i_1}\leq \bar{z}_{i_1}^k-\epsilon_{i_1}. \label{mFPA4a}
\end{equation}

\begin{center}
\begin{minipage}{15cm}
{\footnotesize
\begin{algorithm}[H] {
 \SetKwData{In}{\textbf{Input}} 
   \SetKwData{Out}{\textbf{Output}}
   \SetKwData{Calc}{\textbf{Compute}}
   \SetKwData{Resol}{\textbf{Solve}}
   \SetKwData{Atual}{\textbf{Update}}
   \SetKwData{Const}{\textbf{Build}}
   \SetKwData{Ret}{\textbf{Goback}}
   \SetKwData{Faca}{\textbf{Do}}
   \SetKwData{Par}{\textbf{Stop!}}
   \SetKwData{Cha}{\textbf{Call}}
    \SetKwData{Ins}{\textbf{Insert}}
  
   \SetKwData{Calc}{\textbf{Calculate}}
   \SetKwData{Resol}{\textbf{Solve}}
   \SetKwData{Atual}{\textbf{Update}}
   \SetKwData{Const}{\textbf{Build}}
   \SetKwData{Ret}{\textbf{Go back}}
   \SetKwData{Faca}{\textbf{Do}}
   \SetKwData{Par}{\textbf{Stop!}}
   \SetKwData{Cha}{\textbf{Call}}
   \SetAlgoLined
	\In {$\mathcal{Z}^{*}=\{ \}$, $A, z^{L1}, z^{L2}, z^{I} $,
	          $ \gamma_{i}>0, \zeta_i\in (0,\gamma_{i}), \epsilon_i\in (0,\gamma_i], i=1,2$;  and the permutation  $(i_1,i_2)$}\;
	\Out{$\mathcal{Z}^{*}$ - The Pareto front or an approximation}\;
	\eIf {$z^{L1}=z^{L2}$}
            {\Par The problem has only one nondominated point.
            }
	    {\For{$i=1,2$}
               {
                \If{$i=i_1$}{$w_i=\dfrac{\gamma_i - \zeta_i}{|z^{L1}_{i_1}-z^{L2}_{i_1}|}$}
                \If{$i=i_2$}{$w_i=1$}
                }
            $k=0$\;
	 \If{$DCG$}{\Resol the linear relaxation of subproblem $S^0$ by column generation (Algorithm \ref{AlgGC}) and get the updated $A$ matrix} 
          $\bar{z}^{0}=\left(\sum\limits^{ncol(A)}_{j=1}\bar{x}^{0}_j,\sum\limits^{ncol(A)}_{j=1}\bar{y}^{0}_j \right)$\;
	  \Atual    $\mathcal{Z}^{*} \leftarrow \mathcal{Z}^{*} \cup \{\bar{z}^{0}\}$\; 
	\While{$f_{i_1}(\bar{x}^{k},\bar{y}^{k}) > f^{I}_{i_1}$}
                {
                  \If{$DCG$}{\Resol the linear relaxation of subproblem $S^{k+1}$ by column generation (Algorithm \ref{AlgGC}) and get the updated  $A$ matrix}
	           \Resol $S^{k+1}$ and get $(\bar{x}^{k+1},\bar{y}^{k+1})$ and  
                        $\bar{z}^{k+1}=\left(\sum\limits^{ncol(A)}_{j=1}\bar{x}^{k+1}_j,\sum\limits^{ncol(A)}_{j=1}\bar{y}^{k+1}_j \right)$\;
                  \Atual $\mathcal{Z}^{*} \leftarrow \mathcal{Z}^{*} \cup \{\bar{z}^{k+1}\}$\;
				 $k \leftarrow k+1$
                } 
	 }   
  }
	\caption{ $FPA$ algorithm to solve the  $BOCSP$\label{algFPA}}
\end{algorithm}
}
\end{minipage}
\end{center}

\subsection{The Augmented Weighted Tchebycheff Method} \label{TPA}
The Tchebycheff's method uses a scalarization obtained by minimizing the weighted distance between a feasible solution and the ideal vector with respect to Tchebycheff's norm ($l_{\infty}$). In this study we used a variation of the Tchebycheff Method proposed in \cite{Angelo2018} to overcome the difficulties of calculating the weights and considering a linearization of the objective function. It has been shown (\textit{e.g.} \cite{Angelo2018}) that if the weights are strictly positive, then the solution to the scalarized problem is efficient, moreover, if the solution is unique then it is efficient.

The Augmented Weighted Tchebycheff Method (AWT) ensures that all efficient solutions are found, even those that are in the nonconvex region of the Pareto front, by making an adequate variation of the weights. Considering the lexicographic points, the idea is to normalize the objective functions using  constants $\beta_i, i=1, 2$, calculated via equation \eqref {fatonorm} so that the variation of weights belongs to the interval $[0, 1]$.
    \begin{equation}\label{fatonorm}
    \beta_{i}=\dfrac{1}{|f_{i}^{L2}-f_{i}^{L1}|}, i=1,2.
    \end{equation}
Thus, given a value $\triangle$ (step size) calculated according to  \eqref{passotche} as suggested in \cite{Angelo2018}, at each iteration $w$ is updated to $w = w-\Delta>0$. The initial value for $w$ is defined as $w = 1-\triangle$. If $\triangle=1$ the Pareto front contains only the lexicographic points.
\begin{equation}\label{passotche}
\triangle=\dfrac{1}{f_{2}^{L1}-f_{2}^{L2}}.
\end{equation}

The ISP for the AWT method consists in solving  \eqref{Mtche0}-\eqref{Mtche6} for each given $w$ value. 
\begin{eqnarray}
min &u+\rho\left(\beta_1\left(\sum_{j=1}^{n}x_j-f^{I}_1\right)+\beta_2\left(\sum_{j=1}^{n}
y_j-f^{I}_2\right)\right)\label{Mtche0} \\
s.t. &  (x,y) \in \mathcal{X} \cap \{\eqref{Mtche4} \} \cap \{\eqref{Mtche5}\} \\ 
\quad & \quad u\in \mathbb{R}_{+}.\label{Mtche6}
\end{eqnarray}
\begin{eqnarray}
 & \beta_1 w\left(\sum_{j=1}^{n}x_j-f^{I}_1\right)\leq u \label{Mtche4}\\
 & \beta_2 (1-w)\left(\sum_{j=1}^{n}y_j-f^{I}_2\right)\leq u \label{Mtche5}
\end{eqnarray}

 If the dynamic column generation is employed ($DCG$) the Algorithm \ref{AlgGC} is applied to solve  the subproblem \eqref{Mtche0}-\eqref{Mtche6} for each value of the parameter $w$.\ Otherwise, all the subproblems are solved considering the same set of columns (set $CG_a$ as described in Section \ref{secGCd}), the Static Column Generation method (SCG).

\section{Computational study} \label{comp}
In this section the results of the computational study of  the Lexicographic $\epsilon$-constraint method ($LEC$), the  Frontier Partitioner Algorithm ($FPA$), and the Augmented Weighted Tchebycheff Method ($AWT$)  for the BOCSP are presented.  The algorithms were implemented in the Julia language (Version 0.6.3.1) \cite{Julia} and the mathematical models were implemented using the Modeling Language JuMP \cite{JuMP}. All the subproblems involved were solved using the optimization system CPLEX version 12.6.1 \cite{IBM2014}.  The runs were executed on a computer with an Intel Core i7-3770 processor (3.40GHz), 12 GB of RAM, under a 64-bit Windows Operating System. 

For the pricing subproblems (solved during the column generation algorithm) it was imposed a maximum execution time for CPLEX of 15s and/or a gap of 0.01 and for all the others subproblems (solved during the scalarization's algorithms) it was imposed  a  maximum execution time  of 60s and/or a gap of 0.0001.  Due to these limitations it cannot be ensured that the optimal Pareto fronts were obtained since there is  no guarantee that the  subproblems  were solved to  optimality. Therefore we refer to the set of nondominated points obtained by each run as a \textit{Pareto Front Approximation} (FrA). In addition, a limitation of up to five consecutive iterations of the column generation algorithm with the same dual variable was  defined to  avoid  non-convergence of the algorithm.

In Section \ref{Res1D} we discuss the results for the 1D-BOCSP considering in the FPA Algorithm (Algorithm \ref{algFPA}) the permutation $\{2,1\}$ and $\zeta_i=0.3, i=1,2$. For the 2D-BOCSP, Section \ref{Res2D}, we discuss the results of two variations of the FPA Agorithm. For $c=7$ the results were obtained considering the permutation $\{1,2\}$ and $\zeta_i=0.3, i=1,2$ and for $c=d_{max}$ the results were obtained considering the permutation $\{2,1\}$ and  $\zeta_i=0.1, i=1,2$. These are the configurations that provided the best results.

The metrics $\sigma^1$ (cardinality of the front), $\sigma^2$ (hypervolume), $\sigma^3_{o}$ and $\sigma^3_{c}$ (amplitude of the lexicographic points associated with the minimization of objects and minimization of cycles respectively), $\sigma^{4}$ (total number of subproblems solved), $\sigma^5=\dfrac{\sigma^1}{Total time}$ (number of solutions per execution time) and $\sigma^6=\dfrac{\sigma^1}{\sigma^4}$ (number of points per number of subproblems solved) were used to evaluate the results. The performance profile \cite{Dolan2002} and  boxplot  were also employed in the analyses.

\subsection{The one-dimensional case}\label{Res1D}
The computational study for 1D-BOCSP  was structured in two parts. At first, we show the efficiency of Dynamic Column Generation (DCG) when compared to Static Column Generation (SCG) in the implementation of the scalarizations $LEC, AWT$, and $FPA$. Then we show how the approximation of the Pareto front is affected by the saw capacity value and present a general evaluation of the three scalarizations considering the DCG strategy. 

The 1D instances were generated using CUTGEN1 \cite{GAU1995} considering (v1, v2) = (0.01, 0.2) for Small items (S); (v1, v2) = (0.01,  0.8) for items of small and Medium sizes (M); and (v1, v2) = (0.2, 0.8) for medium and large (Grand)  items (G). We generated three instances of each type (S, M, G) considering $m = 100, L = 10000$, and average demand ($\bar{d}$) equal to  100. In spite of setting $m=100$, CUTGEN1 returned $m<100$ (95 items for the type S, 99 items for the type M and 97 items for the type G) due to duplicated items. For each of these three instances, five new ones were obtained considering subsets of items ($m = 10, 20, 40, 60, 80$). We also considered three values for the saw capacity ($c=4, 7, d_{max}$) resulting in a total of  54 instances, 18 for each  type of item ($S, M$ and $G$). 

For 34 out of the 54 instances tested, the algorithms returned only the lexicographic points ($\sigma^1<=2$). The results shown in Tables \ref{PCE1Dc7}-\ref{conjinstancias}, in  Figures \ref{ResFPA} and \ref{boxplotPCE1D} refer to the instances for each at least one of the algorithms returned points besides the lexicographics ones ($\sigma^1>2$). 

\subsubsection*{Dynamic  versus Static  column generation} \label{dinXest}
Table \ref{PCE1Dc7} presents the results obtained considering the  Static Column Generation (SCG) and the  Dynamic Column Generation (DCG) for the three methods and $c=7$. For each case, the columns display the total number of generated columns (nc), the total computational time (tt), the total number of iterations (it), and the cardinality of the approximations of Pareto front ($\sigma^1$). The results are shown only for the six instances that at least one of the algorithms returned points besides the lexicographics ones ($\sigma^1>2$).
\begin{table}
\small
\begin{center}
\caption{Static (SCG) versus Dynamic (DCG)  column generation for the  $1D-BOCSP$, $c=7$. }
\begin{tabular} {p{1.1cm} | p{1cm}| rrrrrrrrr}
\hline
\multirow{2}{*}{  Method} &\multirow{2}{*}{  id/m}& \multicolumn{2}{c}{nc} & \multicolumn{2}{c}{tt} & \multicolumn{2}{c}{it} & \multicolumn{2}{c}{$\sigma^1$}  \\
&  & SCG & DCG & SCG & DCG & SCG & DCG & SCG& DCG  \\
\hline 
\multirow{6}{*}{LEC } 
& S/40 & 134 & 154 & 3032.72 & 2290.54 & 12 & 13 & 9 & 6\\
& S/60 & 153 & 178 & 1491.21 & 2104.47 & 10 & 10 & 8 & 7\\
& S/80 & 175 & 197 & 926.24 & 928.10 & 5 & 5 &     5 & 5 \\
& S/95 & 207 & 248 & 760.61 & 931.85 & 4 & 5 &     4 & 4\\
& G/60 & 78 & 78 & 23.73 & 23.99 & 3 & 3 &         3 & 3\\
& G/97 & 292 & 306 & 33.54 & 90.44 & 5 & 3 &       5 & 3 \\
\hline
\multirow{6}{*}{FPA } 
& S/40 & 134 & 154 & 960.41 & 761.66 & 12 & 10 &   5 & 5\\
& S/60 & 153 & 178 & 1116.57 & 1093.70 & 15 & 14 & 7 & 8\\
& S/80 & 175 & 197 & 935.29 & 998.42 & 11 & 12 &   5 & 5 \\
& S/95 & 207 & 251 & 825.90 & 697.53 & 10 & 7 &    6 & 4 \\
& G/60 & 78 & 78 & 31.28 & 32.13 & 3 & 3 &         3 & 3 \\
& G/97 & 292 & 306 & 41.13 & 39.09 & 6 & 3 &       5 & 3 \\
\hline
\multirow{6}{*}{AWT} 
& S/40 & 134 & 154 & 1953.60 & 1787.79 & 26 & 26 & 7 & 8\\
& S/60 & 153 & 178 & 2552.89 & 2372.22 & 30 & 35 & 11&  10 \\
& S/80 & 175 & 197 & 1834.99 & 2082.42 & 26 & 30 & 9 &  13\\
& S/95 & 207 & 245 & 2456.16 & 2989.69 & 37 & 45 & 13& 14 \\
& G/60 & 78 & 78 & 24.28 & 25.01 & 1 & 1 &          3& 3 \\
& G/97 & 292 & 306 & 166.49 & 90.15 & 3 & 1 &       5& 3  \\
\hline
\end{tabular}
\label{PCE1Dc7}
\end{center}
\end{table}

For all the instances displayed in Table \ref{PCE1Dc7}, the DCG option resulted in a number of columns greater than or equal to the number of columns generated in the SCG, with the maximum increase in computational time for DCG being $169\%$($id= G, m=97$) in relation to SCG by the $LEC$ method. Considering the $FPA$ and $AWT$ method the maximum increase was $7\%$ ($id=P, m=80$) e $22 \%$ ($id=S, m=95$), respectively. For some instances there was a computational time reduction of up to $45\%$ when considering the dynamic column generation due to the reduction in the number of iterations. As for the cardinality of the FrA ($\sigma^1$) the difference is not significant (\textit{p-value}=$0.75$). 

\begin{figure}
\center
\subfigure[perGCFPA7][Performance Profile of $\sigma^2$  associated with the approximation of the Pareto front obtained by the  $FPA$ method considering DCG and SCG for the $1D-BOCSP$ 
with $c=7$.]{\includegraphics[width=7cm]{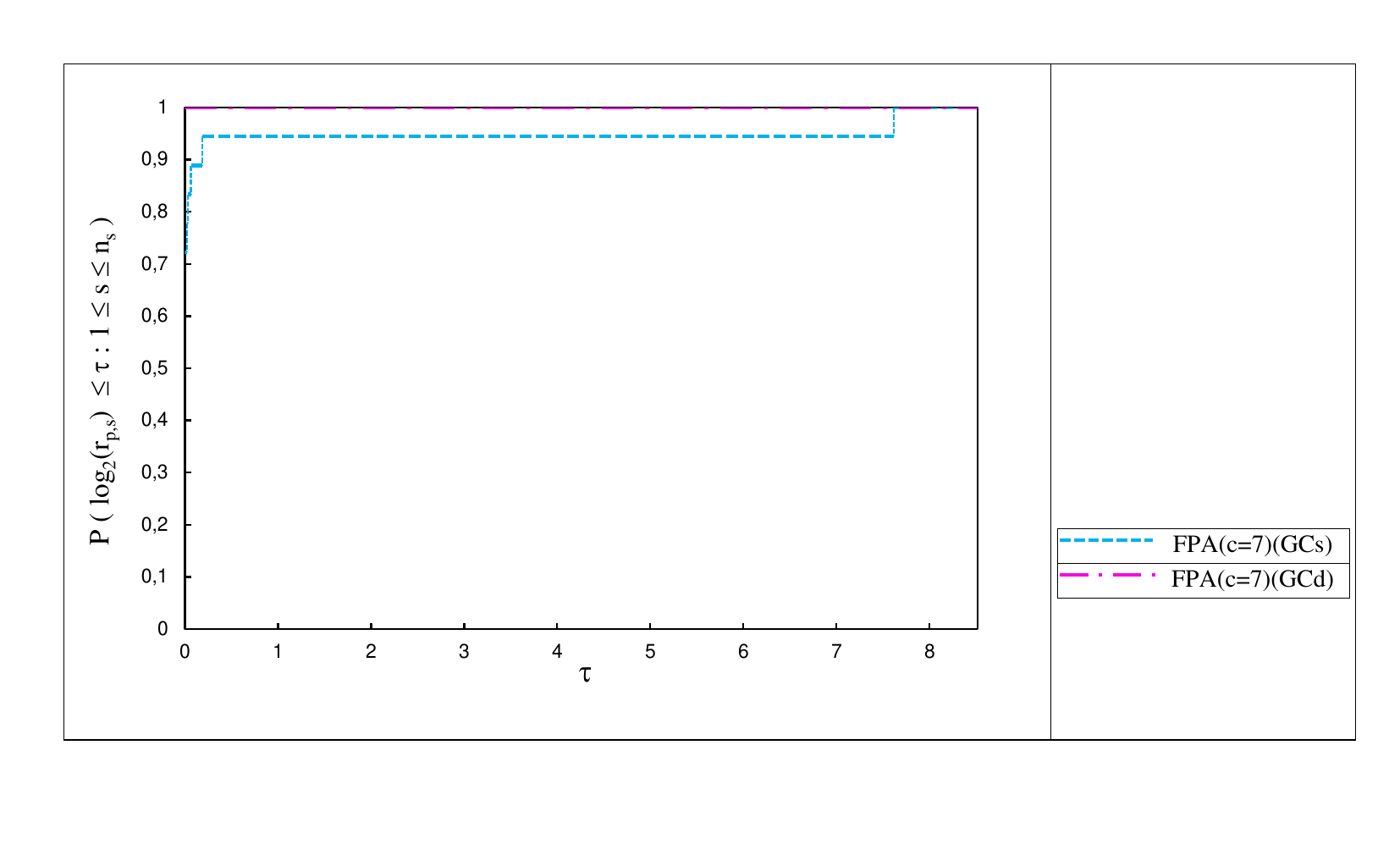}\label{perGCFPA7}}
\qquad
\subfigure[frontPm60][Approximation of the Pareto front obtained by the $FPA$ for the cases DCG and SCG when solving  instance $S/60$, with $c=7$.]{\includegraphics[width=7cm]{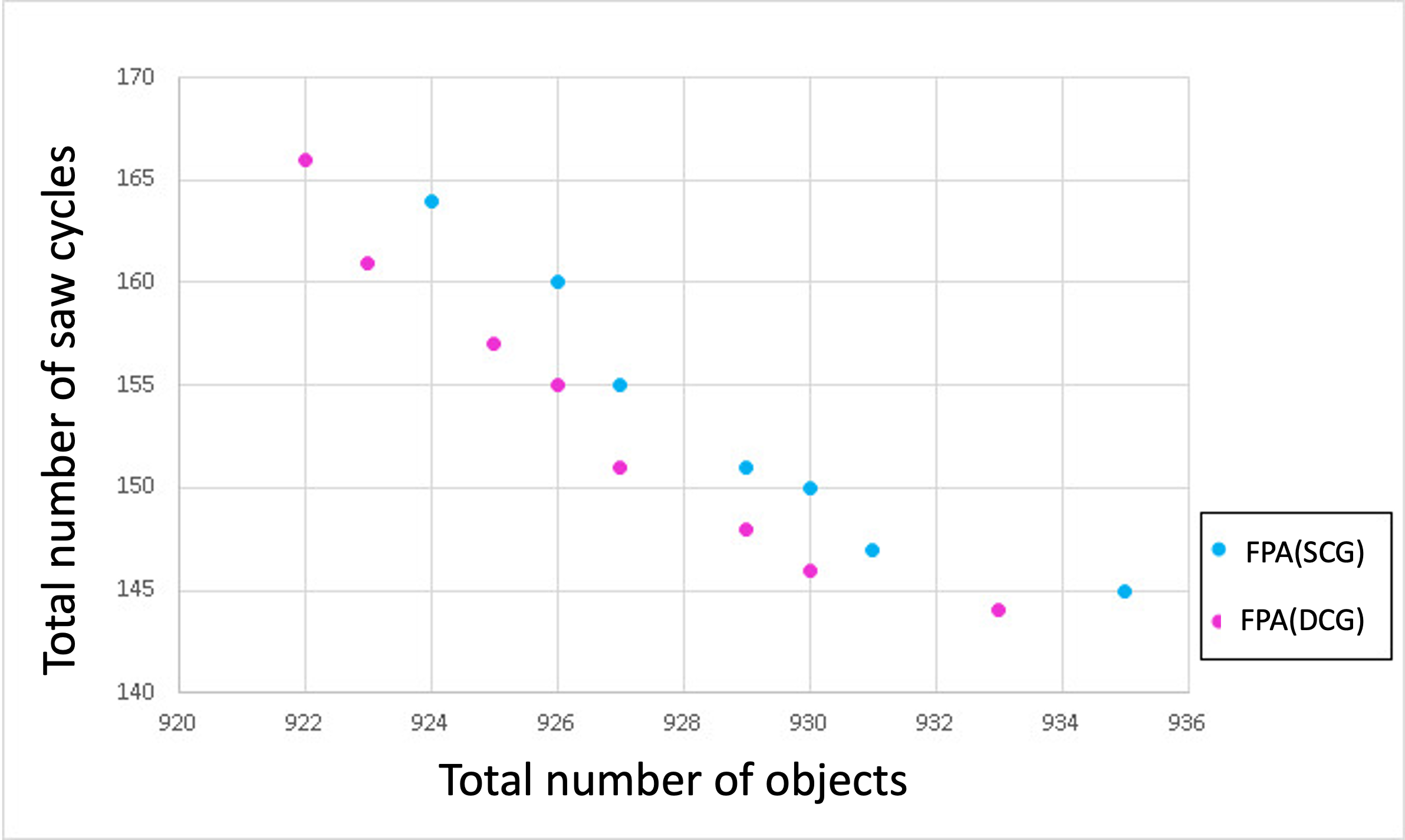}\label{front1DPm60c7GCdxGCs}}
\caption{Dynamic versus Static Column Generation}
\label{ResFPA}
\end{figure}

The quality of the FrA generated by the two strategies of column generation can be better evaluated considering the hypervolume metric ($\sigma^2$). Figure \ref{perGCFPA7} shows the performance profile relative to the results of the FPA algorithm. The performance profile associated with the LEC and AWT algorithms are similar. Employing dynamic column generation provided the best performance for $100\%$ of instances. Figure \ref{front1DPm60c7GCdxGCs} presents the Pareto front approximations obtained applying the $FPA$ algorithm considering the SCG and DCG to solve instance $S/60$ with $c=7$. Note that the front obtained using the DCG is lower and has one more nondominated point, that is, all points dominate the  points obtained using the SCG. This result confirms that the DCG provided a better approximation of the front. Therefore, the Dynamic Column Generation method is employed in the $LEC, FPA$ and $AWT$ algorithms for the remainder of the computational study.

\subsubsection*{Approximation of the Pareto front as the saw capacity increases} \label{card-ciclos}

As mentioned in Section \ref{mod}, the conflict between the minimization of saw cycles and the minimization of objects depends on the scenario. Besides the demand, the saw capacity might also affect the conflict between these two objectives. We study this aspect of the problem testing three different values for the saw capacity ($c = 4, 7, d_{max}$). Table \ref{compa2a} shows the values of the metrics cardinality ($\sigma^1$), hypervolume ($\sigma^2$) and total number of subproblems solved ($\sigma^4$) associated with the FrA obtained by the scalarizations LEC, FPA and AWT (all three with DCG). Only the results for the instances with $\sigma^1 >2$ are shown. We observe that the greater the saw's capacity, the greater is the cardinality of the FrA found by the three methods. The same happens in relation to the hypervolume, and the total number of subproblems solved. For instances of average size ($id=M$) and $m \geq 60$, the methods manage to get an FrA only when $c=d_{max}$. As for $id=G$ there is no FrA for $c=4$. When $c=d_{max}$ the criterion of minimization of cycles is equivalent to minimization of setups, and the results obtained are similar to those in the literature (\textit{e.g.} \cite{Angelo2021}). 
\begin{table}[h!]
\small
\begin{center}
\caption{Metrics $\sigma^1$, $\sigma^2$ and $\sigma^4$ of $LEC$, $FPA$ and $AWT$ with DCG for the  $1D-BOCSP$ as the saw's capacity value changes.}
\begin{tabular}{cc|rrrrrrrrr}
\hline
\multirow{2}{*}{ id/$m$ } & \multirow{2}{*}{ $c$ } & \multicolumn{3}{c}{$\sigma^1$} & \multicolumn{3}{c}{$\sigma^2$} & \multicolumn{3}{c}{$\sigma^4$}\\
& & $LEC$ & $FPA$ & $AWT$ & $LEC$ & $FPA$ & $AWT$ & $LEC$ & $FPA$ & $AWT$  \\
\hline
\multirow{3}{*}{S/40 } 
 &4&5&6&6&489 &491& 503&14&12&22\\
 &7&7&5&8&668 &662 &678&30&14&30\\
 & $d_{max}$ &23&22&18&30141& 30175 &30049&66&32&42\\
 \hline
 \multirow{3}{*}{S/60 } 
&4&5&4&10&650 &683 &780&14&10&28\\
&7&8&8&10&1248 & 1216 &1325&24&18&39\\
& $d_{max}$ &33&36&26&99847& 99845& 99499&90&48&67\\
 \hline
 \multirow{3}{*}{S/80 } 
&4&3&5&10&805 &918& 1061&10&12&33\\
&7&7&5&13&2380&1527 &1814&14&16&34\\
& $d_{max}$ &40&39&30&106019 &106047& 105691&96&52&67\\
\hline
\multirow{3}{*}{S/95 } 
&4&3&4&9&364&361&549&10&10&33\\
&7&5&4&14&2047& 1757 &2103&14&11&49\\
& $d_{max}$ &42&41&32&205682 &205697 &205217&108&60&86\\
\hline
M/60 & $d_{max}$ &9&9&7&5111 &5111 &5067&22&13&10\\
\hline
M/80 & $d_{max}$ &20&21&15&16534 &16587 &16365&44&30&23\\
\hline
M/99 & $d_{max}$ &27&28&20&39782 &39864 &39320&72&43&37\\
\hline
\multirow{2}{*}{G/60 } 
&7&3&3&3&217 &217& 216&10&7&5\\
& $d_{max}$ &5&5&4&1531 &1531 &1530&14&9&7\\
\hline
G/80 & $d_{max}$ &17&17&12&9492 &9493 &9437&38&23&19\\
\hline
 \multirow{2}{*}{G/97 } 
&7&3&3&3&152&152&152&10&7&5\\
& $d_{max}$ &22&24&17&18403 &18492 &18311&56&34&27\\
\hline
\end{tabular}\\
\label{compa2a}
\end{center}
\end{table}

Table \ref{compa1} presents the total number of columns generated  (nc), the total computational time in seconds (tt) and the total number of iterations (it) for each method. The total number of columns generated by all the three methods is the same for 15 out of the 20 instances. The difference is smaller than $4\%$ for the four small instances (id=S) with $c=4$ and one small instance with $c=7$. However the total computational time and total number of iterations are quite different. The total computational time  spent by the AWT method can be up to 6.4 times greater than the $FPA$ method (instance G/80 c=$d_{max}$) and 3.2 times greater than the LEC method (instance S/95 c=7).  And the total number of iterations that the $AWT$ method uses to obtain an FrA can be up to 9.6 times greater than $LEC$ (instance S/95 and $c=4$) and 6.4 times greater than  $FPA$ (instance S/97 and $c=7$). It is interesting to note that, in spite of the greater computational time and number of iterations, the cardinality of the FrA given by the $AWT$ is not always greater than the others, this is the case, for example  of the instance G/80 with c=$d_{max}$ as shown in Table \ref{compa2a}. 

\begin{table}[h!]
\small
\begin{center}
\caption{Computational effort of $LEC$, $FPA$ and $AWT$ with DCG for the $1D-BOCSP$.}
\begin{tabular}{cc|rrrrrrrrrr}
\hline
\multirow{2}{*}{ id/$m$ } & \multirow{2}{*}{ $c$ } & \multicolumn{3}{c}{nc} & \multicolumn{3}{c}{tt} & \multicolumn{3}{c}{it} \\
 & & $LEC$ & $FPA$ & $AWT$ & $LEC$ & $FPA$ & $AWT$ & $LEC$ & $FPA$ & $AWT$  \\
\hline \hline
S/40 & 4 & 134 & 136 & 134 & 641.14 & 643.70 & 1245.41 & 5 & 8 & 18 \\
S/40 & 7 & 154 & 154 & 154 & 2290.54 & 761.66 & 1787.79 & 13 & 10 & 26 \\
S/40 & $d_{max}$ & 148 & 148 & 148 & 3164.51 & 1381.61 & 2117.83 & 31 & 28 & 38 \\
\hline
S/60 & 4 & 170 & 174 & 170 & 927.92 & 634.42 & 2271.91 & 5 & 6 & 24 \\
S/60 & 7 & 178 & 178 & 178 & 2104.47 & 1093.70 & 2372.22 & 10 & 14 & 35 \\
S/60 & $d_{max}$ & 178 & 178 & 178 & 3339.27 & 1906.62 & 3401.67 & 43 & 44 & 63 \\
\hline
S/80 & 4 & 191 & 198 & 191 & 686.91 & 754.79 & 2023.54 & 3 & 8 & 29 \\
S/80 & 7 & 197 & 197 & 197 & 928.10 & 998.42 & 2082.42 & 5 & 12 & 30 \\
S/80 & $d_{max}$ & 191 & 191 & 191 & 2969.49 & 1849.75 & 3020.45 & 46 & 48 & 63 \\
\hline
S/95 & 4 & 225 & 226 & 225 & 688.32 & 636.97 & 2022.69 & 3 & 6 & 29 \\
S/95 & 7 & 248 & 251 & 245 & 931.85 & 697.53 & 2989.69 & 5 & 7 & 45 \\
S/95 & $d_{max}$ & 240 & 240 & 240 & 3723.53 & 2200.61 & 4078.03 & 52 & 56 & 82 \\
\hline
M/60 & $d_{max}$ & 175 & 175 & 175 & 46.88 & 40.51 & 105.84 & 9 & 9 & 6 \\
\hline
M/80 & $d_{max}$ & 294 & 294 & 294 & 1595.46 & 1133.96 & 1153.62 & 20 & 26 & 19 \\
\hline
M/99 & $d_{max}$ & 379 & 379 & 379 & 3736.20 & 2268.57 & 2007.88 & 34 & 39 & 33 \\
\hline
G/60 & 7 & 78 & 78 & 78 & 23.99 & 32.13 & 25.01 & 3 & 3 & 1 \\
G/60 & $d_{max}$ & 78 & 78 & 78 & 24.27 & 31.63 & 24.95 & 5 & 5 & 3 \\
\hline
G/80 & $d_{max}$ & 328 & 328 & 328 & 373.85 & 104.99 & 672.23 & 17 & 19 & 15 \\
\hline
G/97 & 7 & 306 & 306 & 306 & 90.44 & 39.09 & 90.15 & 3 & 3 & 1 \\
G/97 & $d_{max}$ & 431 & 431 & 431 & 2605.68 & 1486.87 & 1426.86 & 26 & 30 & 23 \\
\hline
\end{tabular}\\
\vspace{0.2cm}
\label{compa1}
\end{center}
\end{table}

\begin{figure}[h!]
\center
\subfigure[Sigma5][Boxplot relative to $\sigma^5$ for the $1D-BOCSP$.]{\includegraphics[width=7cm]{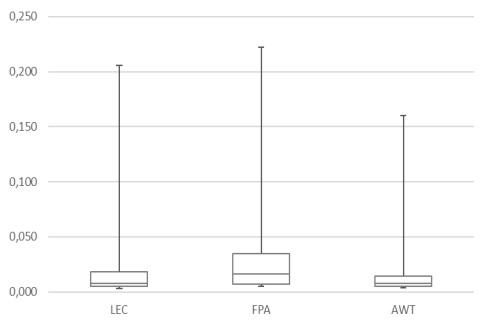}\label{sig5-1D}}
\qquad
\subfigure[sigma6][Boxplot relative to $\sigma^6$ for the $1D-BOCSP$.]{\includegraphics[width=7cm]{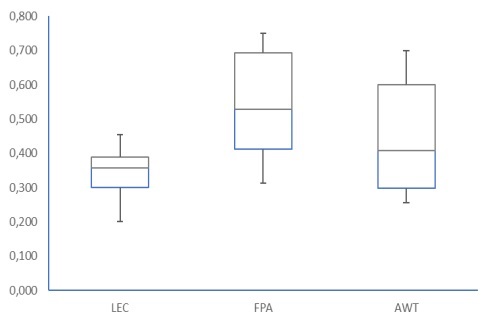}\label{sig6-1D}}
\caption{Computational effort to obtain an FrA for the 1D-BOCSP relative to $\sigma^5$ and $\sigma^6$.}
\label{boxplotPCE1D}
\end{figure}

As for the number of nondominated points by execution time ($\sigma^5$) and  by subproblem ($\sigma^6$) the results  presented in Table \ref{compa5} show that the $FPA$ method obtained the greatest values in more than 75\% of the instances. However its variability  is greater than the other two methods considering $\sigma^5$ and similar to AWT considering $\sigma^6$  as shown in Figure \ref{boxplotPCE1D}. 

\begin{figure}
\centering
\includegraphics[width=10cm,height=6cm]{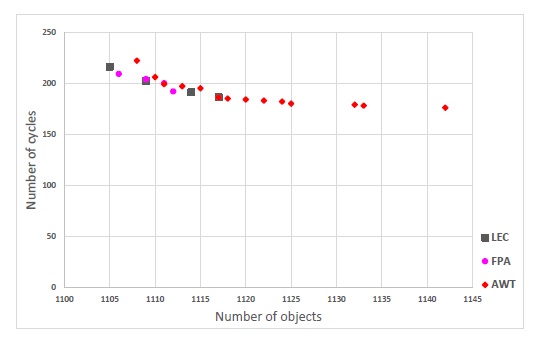}
\caption{Approximation of the Pareto Front for the instance S/95 of the 1D-BOCSP with  $c=7$ (with DCG).}

\label{frontmetc7}
\end{figure}

Figure \ref{frontmetc7} shows the FrA obtained for the instance $S/100$ with $c=7$ by each of the three methods. The $FPA$ method generated fewer nondominated points, however they dominate the ones obtained by the two other methods. The FrA obtained by the  $AWT$ method has points spread along the rectangle that contains the Pareto front, thus obtaining a greater hypervolume ($\sigma^2$). The FrAs shown in Figure \ref{frontmetc7} indicate that a combination of these three methods might provide a more interesting FrA. This observation is further explored in Section \ref{Res3}.

\begin{table}[h!]
\small
\begin{center}
\caption{Metrics $\sigma^5$ e $\sigma^6$ of $LEC$, $FPA$ e $AWT$ for the $1D-BOCSP$ and  $c=4, 7, d_{max}$.}
\begin{tabular}{cc|rrrrrr}
\hline
\multirow{2}{*}{ id/$m$ } & \multirow{2}{*}{ $c$ } & \multicolumn{3}{c}{$\sigma^5$} & \multicolumn{3}{c}{$\sigma^6$} \\
 & & $LEC$ & $FPA$ & $AWT$ & $LEC$ & $FPA$ & $AWT$ \\
\hline \hline
S/40 & 4         & 0.008 & 0.009 & 0.005 & 0.357 & 0.500 & 0.273\\
S/40 & 7         & 0.003 & 0.007 & 0.004 & 0.200 & 0.357 & 0.267 \\
S/40 & $d_{max}$ & 0.007 & 0.016 & 0.008 & 0.348 & 0.687 & 0.429 \\
\hline
S/60 & 4         & 0.005 & 0.006 & 0.004 & 0.357 & 0.400 & 0.357 \\
S/60 & 7         & 0.003 & 0.007 & 0.004 & 0.292 & 0.444 & 0.256 \\
S/60 & $d_{max}$ & 0.010 & 0.019 & 0.008 & 0.367 & 0.750 & 0.388 \\
\hline
S/80 & 4         & 0.004 & 0.008 & 0.005 & 0.300 & 0.417 & 0.303 \\
S/80 & 7         & 0.005 & 0.005 & 0.006 & 0.357 & 0.312 & 0.382 \\
S/80 & $d_{max}$ & 0.013 & 0.021 & 0.010 & 0.406 & 0.750 & 0.448 \\
\hline
S/95 & 4         & 0.004 & 0.006 & 0.004 & 0.300 & 0.400 & 0.273 \\
S/95 & 7         & 0.004 & 0.006 & 0.005 & 0.286 & 0.364 & 0.286 \\
S/95 & $d_{max}$ & 0.011 & 0.019 & 0.008 & 0.389 & 0.683 & 0.372 \\
\hline
M/60 & $d_{max}$ & 0.192 & 0.222 & 0.066 & 0.409 & 0.692 & 0.700 \\
\hline
M/80 & $d_{max}$ & 0.012 & 0.018 & 0.013 & 0.454 & 0.700 & 0.652 \\
\hline
M/99 & $d_{max}$ & 0.007 & 0.012 & 0.010 & 0.375 & 0.651 & 0.540 \\
\hline
G/60 & 7         & 0.125 & 0.093 & 0.120 & 0.300 & 0.428 & 0.600 \\
G/60 & $d_{max}$ & 0.206 & 0.158 & 0.160 & 0.357 & 0.556 & 0.571 \\
\hline
G/80 & $d_{max}$ & 0.045 & 0.162 & 0.018 & 0.447 & 0.739 & 0.632 \\
\hline
G/97 & 7         & 0.033 & 0.077 & 0.033 & 0.300 & 0.428 & 0.600 \\
G/97 & $d_{max}$ & 0.008 & 0.016 & 0.012 & 0.393 & 0.706 & 0.630 \\
\hline
\end{tabular}\\
\label{compa5}
\end{center}
\end{table}

\subsection{Results for the two-dimensional case}\label{Res2D}
To define  the set of instances used in the computational study for the 2D-BOCSP, preliminary studies were conducted considering several factors, including the shape of the items, the number of different types of items, and the value of the saw capacity.  In the following it is  reported the computational results of applying the methods LEC, FPA and AWT with dynamic column  generation for c=7 and $c=d_{max}$.  The instances described in Table \ref{conjinstancias} were  generated using the 2DCPackGen \cite{Silva2014} based on the parameters used in \cite{Mateus2020} and the real shape of the items of a small scale furniture industry \cite{Figueiredo}. The instances are named  $id/m$ with  id (number) representing the shape of the item: id= 1 (small and square), id=3 (medium length and narrow), id=6 (large and square), id=11 (medium size and square), and id=14 (long and narrow or long and tall).  The object was generated using id=10 (long and average width). Altogeter 40 instances were generated considering  five values for id, four values for $m$ and two values for the saw capacity (c).
  \begin{table}
  \begin{center}
  \caption{Instances of the $2D-BOCSP$.}
  \begin{tabular}{c| c }
 \hline
Parameters& values \\
  \hline
 \hline
id(items)   &  1, 3, 6, 11 and 14 \\
 \hline
 $m$  &  20, 30, 50 and 100\\
 \hline
 $c$ &  7 and $d_{max}$ \\
 \hline
 id(object)  &  10 \\
 \hline
 $L$ e $W$   & [100,200]\\
 \hline
 $l$ e $w$  &  [25,100]\\
 \hline
 $d$  & [10,200] \\
 \hline
 \end{tabular}\\
 \vspace{0.4cm}
 \label{conjinstancias}
 \end{center}
 \end{table}

\begin{figure}[h!]
\center
\subfigure[Sigma5][Boxplot relative to $\sigma^5$ for the $2D-BOCSP$.]{\includegraphics[width=7cm]{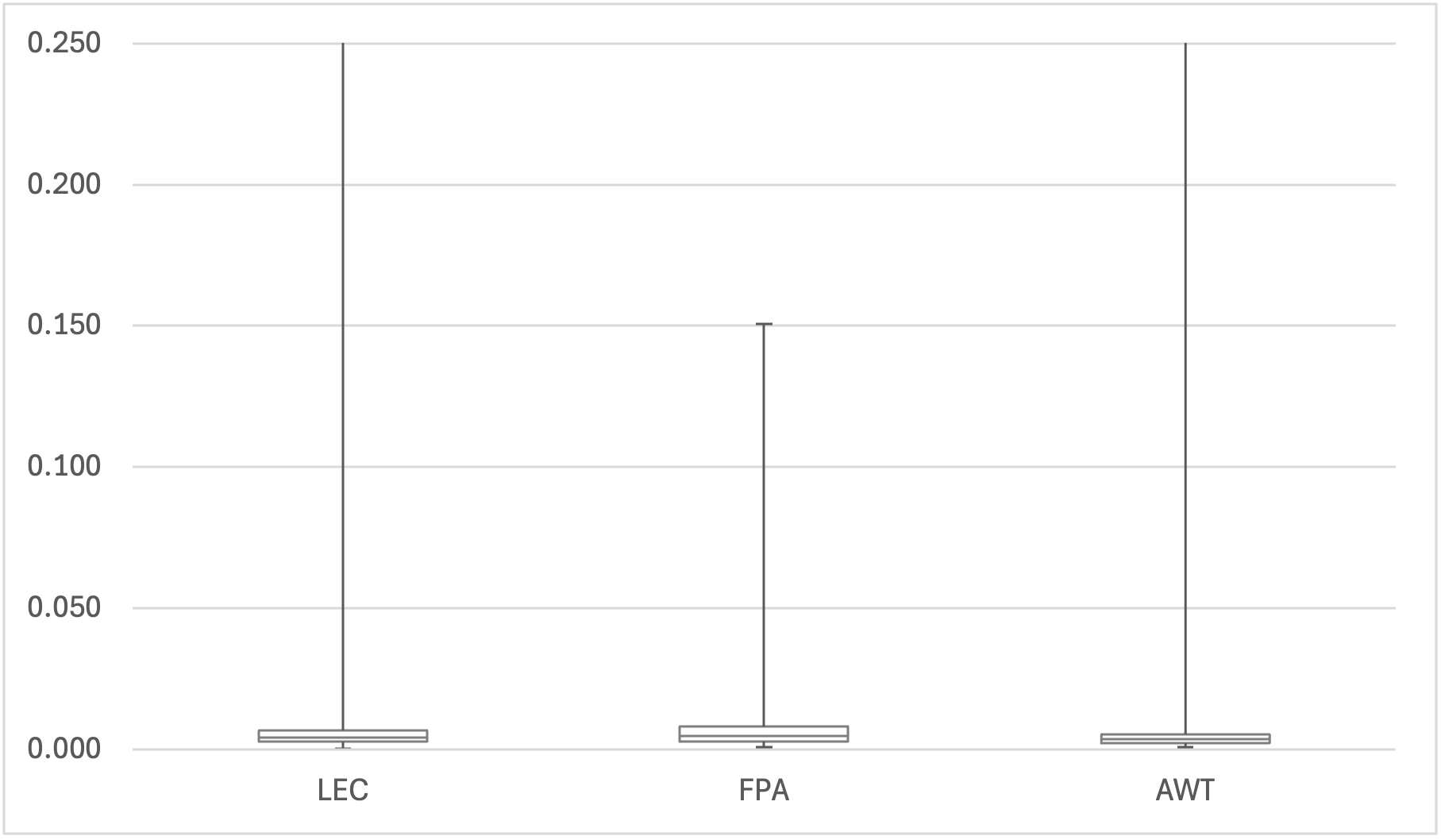}\label{sig5-2D}}
\qquad
\subfigure[sigma6][Boxplot relative to $\sigma^6$ for the $2D-BOCSP$.]{\includegraphics[width=7cm]{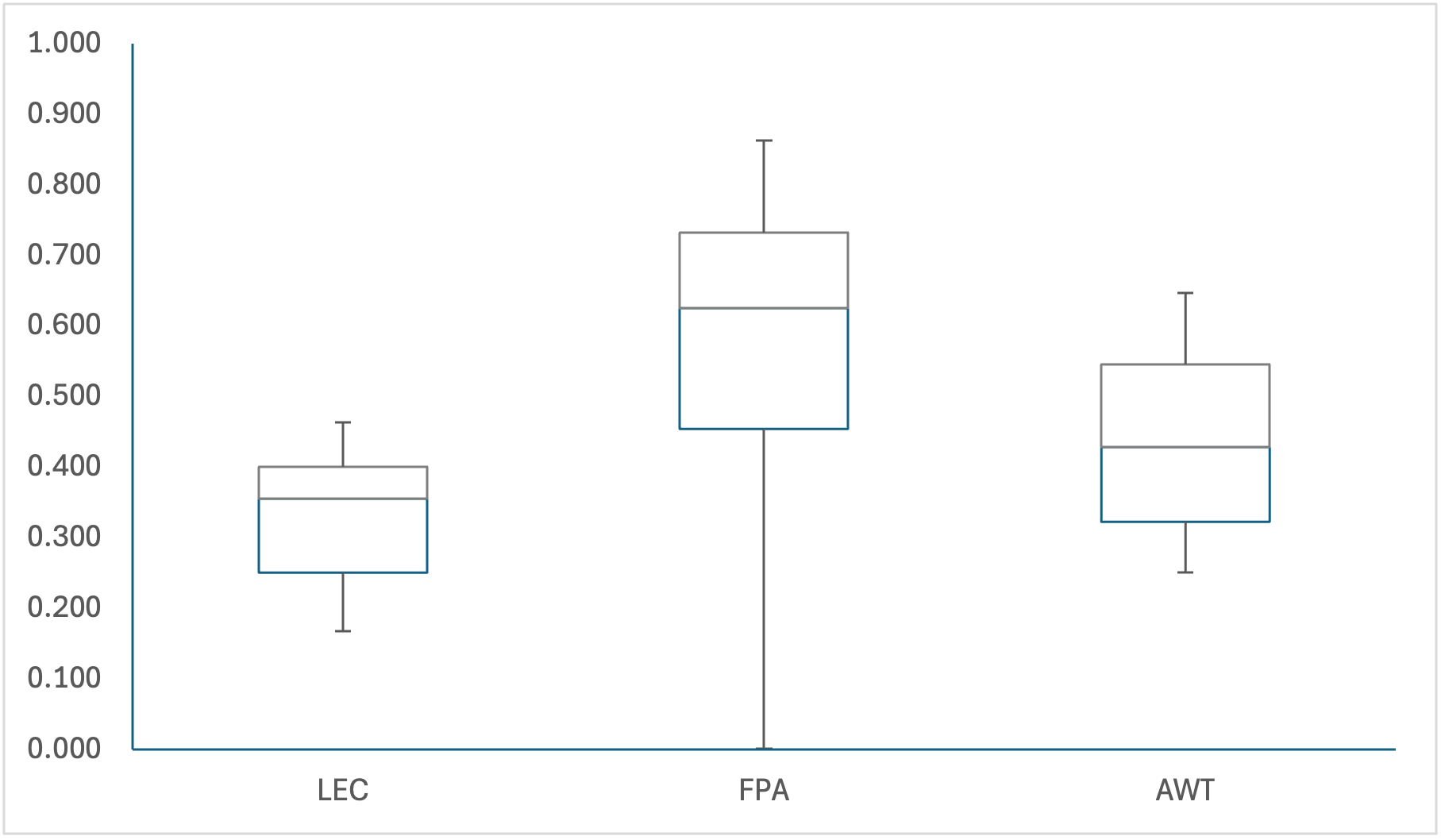}\label{sig6-2D}}
\caption{Computational effort to obtain an FrA for the 2D-BOCSP associated with $\sigma^5$ and $\sigma^6$.}
\label{boxplotPCE2D}
\end{figure}

\begin{figure}[h!]
\centering
\includegraphics[width=10cm,height=6cm]{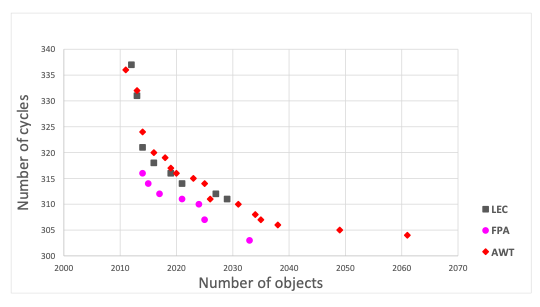}\\
\caption{Approximation of the Pareto front given by the methods  $FPA$, $LEC$ and  $AWT$ with DCG for the instance 11/100 of the 2D-BOCSP with  $c=7$.}
\label{frontmet2Dc7}
\end{figure}

For 9 out of the 40 instances tested, the algorithms returned only the lexicographic points ($\sigma^1\leq 2$). Three of these instances have id=6 (items that are large and square) and the saw capacity equals to  7. The results shown in Tables \ref{sigma1-sigma4-2D-cap7dmax}-\ref{compmetcap7dmax2}, in Figures \ref{boxplotPCE2D} and \ref{frontmet2Dc7}  refer to the instances for which at least one of the algorithms returned points besides the lexicographics ones ($\sigma^1>2$). 

Table \ref{sigma1-sigma4-2D-cap7dmax} shows the values of the metrics $\sigma^1$, $\sigma^2$ and $\sigma^4$ associated
to the FrA obtained by the scalarizations $LEC$, $FPA$ and $AWT$ (all three with DCG).  As in the one-dimensional case, the cardinality of the FrA found by the three methods 
is greater when $c=d_{max}$. The same happens in relation to the hypervolume, and the total number of subproblems solved. The greatest hypervolume for the saw capacity  equal to $d_{max}$ is associated with the  method $FPA$ (instance  $14/100, \sigma^2=80259)$, and in this case it is associated with the smallest number of subproblems ($\sigma^4=67$).  The FrA with the greater cardinality was obtained by the method $LEC$ for the instance $11/100$ ($\sigma^1$=53),  it provided 36\% more points than the method AWT  and 11\% than the method FPA. For this instance the associated hypervolume is also greater (25\% and 27\%) than the ones associated with the methods $FPA$ and $AWT$ respectively). However, the associated total number of subproblems ($\sigma^4$) is 75\% and 100\% greater than the values  associated with the methods $AWT$ and $FPA$ respectively. 

For $c=7$, the greatest difference in $\sigma^1$ was seen for instance $11/100$, where the method $AWT$ generated an FrA with a number of points 70\% greater   than the other two methods ($\sigma^1$ = 16), although, the associated hypervolume was smaller than the one associated with the FPA algorithm that returned only 7 nondominated points with a total number of subproblems ($\sigma^4$) 57\% smaller and  $\sigma^2=2317$, the greater hypervolume among all the three FrA obtained for this instance. This indicates that despite the fact that an FrA has fewer points, it can be closer to the optimal Pareto front. 
 \begin{table}
\small
\begin{center}
\caption{Metrics associated with the FrA obtained by the  methods $LEC$, $FPA$ and  $AWT$ for the  $(PCEB_{c}^{L}-2D)$ as the saw’s capacity value changes (with DCG).}
\begin{tabular}{rr|ccccccccc}
\hline
\multirow{2}{*}{ id/m}&\multirow{2}{*}{ c} & \multicolumn{3}{c}{$\sigma^1$} & \multicolumn{3}{c}{$\sigma^2$} & \multicolumn{3}{c}{$\sigma^4$} \\
&& $LEC$ & $FPA$ & $AWT$ & $LEC$ & $FPA$ & $AWT$ & $LEC$ & $FPA$ & $AWT$ \\
\hline \hline
1/20&7         & 6 & 4 & 5 & 491 &   303 &  315 & 16 & 9 & 14 \\
1/20&$d_{max}$ & 13 & 12   & 10 & 5716 &  5720  & 5664 & 32 & 17 & 22 \\
\hline
1/30&7         & 6 & 7 & 5 &532 &   387 &  394 & 20 & 11 & 18 \\
1/30&$d_{max}$ & 19 & 15 & 14 & 9969  & 10017  &  9916& 48 & 23 & 35 \\
\hline
1/50&7 & 10 & 8 & 10 & 804 &  720 &  748 & 22 & 13 & 19 \\
1/50&$d_{max}$ & 25 & 25 & 19 & 14114 &  14114 &  13887 & 54 & 29 & 30 \\
\hline
3/20&7 & 6 & 4 & 4 & 507 &  325 &  322 & 16 & 8 & 13 \\
3/20&$d_{max}$ & 14 & 14 & 9 &  5665 &   5692 &  5612& 32 & 18 & 19 \\
\hline
3/30&7 & 7 & 4 & 4 & 420 &  398 &  399 & 22 & 10 & 14 \\
3/30&$d_{max}$ & 19 & 15 & 14 & 23634 &  13207 &  12605 & 52 & 24 & 35 \\
\hline
3/50&7 & 6 & 6 & 7 & 564 &  573  & 575 & 18 & 11 & 26 \\
3/50&$d_{max}$ & 27 & 30 & 23 & 32437 &  32425 &  32201 & 76 & 41 & 49 \\
\hline
6/20&$d_{max}$ & 8 & 8 & 6 & 2500 &   2508 &  2493 & 20 & 12 & 10 \\
\hline
6/30&$d_{max}$ & 14 & 14 & 10 & 5485 &  5485 &  5449 & 32 & 18 & 17 \\
\hline
6/50&$d_{max}$ & 14 & 14 & 11 & 6455 &  6476  & 6389 & 32 & 18 & 17 \\
\hline
6/100&7 & 5 & 5 & 4 &  339  &   346  & 358 & 16 & 11 & 14 \\
6/100&$d_{max}$ & 38 & 39 & 26 & 46004 &  46180 &  45808 & 92 & 48 & 49 \\
\hline
11/20&$d_{max}$ & 12 & 12 & 10 &3660 &  3620 &  3560 & 30 & 16 & 17 \\
\hline
11/30&7 & 5 & 4 & 5 &318  & 248 &  259 & 16 & 10 & 13 \\
11/30&$d_{max}$ & 19 & 15 & 12 &9190 &   9414 &  9099 & 50 & 23 & 22 \\
\hline
11/50&7 & 7 & 7 & 10 &1521 &  969 &   983 & 26 & 13 & 36 \\
11/50&$d_{max}$ & 29 & 22 & 18 & 28281 &  28303 &   28183 & 76 & 40 & 43 \\
\hline
11/100&7 & 9 & 7 & 16 &  2029 &  2317 &  2160 & 28 & 15 & 35 \\
11/100&$d_{max}$ & 53 & 47 & 34 & 78296 &  62463 &  61663 & 110 & 55 & 63 \\
\hline
14/20&$d_{max}$ & 10 & 10 & 7 &  3892 &  3926 &  3863 & 24 & 14 & 11 \\
\hline
14/30&7 & 2 & 3 & 2 & 140  & 141 &  142 & 8 & 7 & 5 \\
14/30&$d_{max}$ & 18 & 17 & 12 & 10536 &  10533  & 10420 & 40 & 23 & 22 \\
\hline
14/50&7 & 2 & 2 & 3 &164 &  164  & 165 & 8 & 6 & 7 \\
14/50&$d_{max}$ & 21 & 17 & 13 & 12906 &  12104 &  12094 & 54 & 26 & 22 \\
\hline
14/100&7 & 4 & 10 & 10 & 1318 &   1535 &  1617& 10 & 16 & 31 \\
14/100&$d_{max}$ & 36 & 35 & 29 & 79284 &  80259  & 79668 & 132 & 67 & 77 \\
\hline
\end{tabular}\\
\vspace{0.2cm}
\label{sigma1-sigma4-2D-cap7dmax}
\end{center}
\end{table}
Table \ref{compmetcap7dmax} shows the total number of columns generated (nc), the total computational time in seconds (tt) and the total number of iterations (it) for each method.  Considering the instances with $c=7$, the $FPA$ method generated more columns than the other two methods for $62\%$ of the instances and for the others instances  the total number of columns generated by the three methods is the same. Additionally, the FPA method obtains an FrA  with a smaller or equal number of iterations for $77\%$ of instances. The  $AWT$ method  used an average number of iterations twice as greater as the other two methods to obtain an FrA. 

In the case of $c=dmax$, the method FPA uses less iterations than the method AWT for 11 out of the 18 instances  (in average 14\% less iterations). The  number of iterations of the method LEC is the same as the method FPA for eight instances, and it uses 26\% more iterations than the method FPA for the other 10 instances.
 
\begin{table}
\small
\begin{center}
\caption{Comparison of the computational effort of LEC, FPA and AWT with
DCG for the BOCSP-$2D$.}
\begin{tabular}{rr|ccccccccc}
\hline
\multirow{2}{*}{ id/$m$ }&\multirow{2}{*}{c} & \multicolumn{3}{c}{nc} & \multicolumn{3}{c}{tt($s$)} & \multicolumn{3}{c}{it} \\
    && $LEC$ & $FPA$  & $AWT$ & $LEC$ & $FPA$ & $AWT$ & $LEC$ & $FPA$  & $AWT$ \\
\hline\hline
1/20&7           & 141 & 144 & 141 & 1446.65 & 1106.14 & 2574.89 & 6 & 5 & 10 \\
1/20&$d_{max}$   & 133 & 135 &133      & 1600.42  & 1209.73   & 1406.30 &14    & 13  &18 \\   
\hline
1/30&7& 233 & 238 & 233 & 3508.47 & 3000.62 & 4283.53 & 8 & 7 & 14 \\
1/30&$d_{max}$    & 226 & 230 &226      & 4784.45  & 3642.32   & 4483.54 &22   & 19   &31  \\  
\hline
1/50&7& 151 & 151 & 151 & 2965.84 & 2422.74 & 3987.24 & 9 & 9 & 15 \\
1/50&$d_{max}$    & 151 & 151 &149      & 3256.71  & 3303.01   & 3946.45 &27    & 25  &26  \\  
\hline
3/20&7 & 103 & 109 & 103 & 760.65 & 492.73 & 1347.98 & 6 & 4 & 9 \\
3/20&$d_{max}$    & 109 & 111 &109      &  887.93  & 863.49    & 611.78 &14    & 14  &15  \\   
\hline
3/30&7 & 198 & 203 & 189 & 2202.39 & 1605.23 & 2339.90 & 9 & 6 & 10 \\
3/30&$d_{max}$    & 223 & 217 &202      & 4026.89  & 2581.10   & 3524.71 &45   & 20   &31  \\  
\hline
3/50&7 & 260 & 282 & 260 & 6868.03 & 6369.71 & 8174.31 & 7 & 7 & 22 \\
3/50&$d_{max}$    & 225 & 225 &224      & 10448.87 & 10927.78  & 6679.63 &180    & 37  &45  \\ 
\hline
6/20&$d_{max}$    & 77 & 77   &77       &  25.79   & 53.06     & 18.29 &8     & 8  &6  \\      
\hline
6/30&$d_{max}$    & 128& 128  &128      & 242.92   & 193.26    & 383.60 &14    & 14  &13  \\   
\hline
6/50&$d_{max}$    & 185& 185  &185      & 1776.18  & 1550.66   & 1587.63 &14    & 14  &13  \\  
\hline
6/100&7 & 356 & 356 & 356 & 1244.33 & 1061.55 & 1270.89 & 6 & 7 & 10 \\
6/100&$d_{max}$   & 364& 365  &364      & 7156.91  & 9345.99   & 5163.98 &44     & 44 &45 \\   
\hline
11/20&$d_{max}$   &156 & 155 &151       & 1024.30  & 879.05    & 761.19 &18     & 12 &13  \\   
\hline                                                                              
11/30&7 & 209 & 228 & 204 & 1214.67 & 1223.24 & 1723.26 & 6 & 6 & 9 \\
11/30&$d_{max}$   &242 & 257 &239       & 3431.73  & 2422.03   & 3121.36 &26     & 19 &18 \\   
\hline                                                                              
11/50&7 & 366 & 366 & 366 & 3088.93 & 2285.46 & 3962.55 & 11 & 9 & 32 \\
11/50&$d_{max}$   &334 & 334 &334       & 5865.09  & 3988.42   & 4237.22 &36    & 36  &39  \\  
\hline                                                                              
11/100&7 & 356 & 357 & 356 & 8717.08 & 7716.79 & 11216.70 & 12 & 11 & 31 \\
11/100&$d_{max}$  &337 & 337 &336       & 12722.87 & 10210.19  & 10714.20 &54     & 51 &59 \\  
\hline                                                                              
14/20&$d_{max}$   &128  & 130 &128      & 398.52   & 512.76    & 306.49 &10  & 10    &7 \\     
\hline                                                                              
14/30&7 & 215 & 215 & 215 & 528.47 & 529.82 & 410.06 & 2 & 3 & 1 \\
14/30&$d_{max}$   &240  & 239 &239      & 3695.03  & 3033.94   & 2915.60 &19  & 19    &18 \\   
\hline                                                                              
14/50&7 & 249 & 249 & 249 & 542.20 & 491.87 & 590.04 & 2 & 2 & 3 \\
14/50&$d_{max}$   &257  & 257 &257      & 2976.94  & 1628.78   & 1420.97 &25  & 22    &18 \\   
\hline                                                                              
14/100&7 & 838 & 906 & 838 & 11301.48 & 14308.69 & 13361.87 & 3 & 12 & 27 \\
14/100&$d_{max}$  &661  & 661 &661      & 16451.12 & 12796.03  & 13161.80 &64  & 63    &73 \\   
\hline                                                                              
\hline
\end{tabular}\\
\vspace{0.2cm}
\label{compmetcap7dmax}
\end{center}
\end{table}

\begin{table}
\small
\begin{center}
\caption{Metrics $\sigma^5$ e $\sigma^6$ associated with the FrA given by the methods $LEC$, $FPA$ and  $AWT$ for the $(PCEB_{c}^{L}-2D)$.}
\begin{tabular}{rr|rrrrrr}
\hline
\multirow{2}{*}{id/$m$} &\multirow{2}{*}{c}  & \multicolumn{3}{c}{$\sigma^5$} & \multicolumn{3}{c}{$\sigma^6$} \\
&& $LEC$ & $FPA$ & $AWT$ & $LEC$ & $FPA$ & $AWT$ \\
\hline
1/20&7 & 0.0028 & 0.0036 & 0.0020 & 0.2500 & 0.4444 & 0.3571 \\
1/20&$d_{max}$   &0.0075  & 0.0099 &0.0071  &0.3750  &0.7059  & 0.4545 \\  
\hline
1/30&7 & 0.0011 & 0.0023 & 0.0012 & 0.2000 & 0.6364 & 0.2778 \\
1/30&$d_{max}$   &0.0040  & 0.0041 &0.0031  &0.3958  &0.6522  & 0.4000 \\  
\hline
1/50&7 & 0.0027 & 0.0033 & 0.0025 & 0.3636 & 0.6154 & 0.5263 \\
1/50&$d_{max}$   &0.0077  & 0.0076 &0.0048  &0.4630  &0.8621  & 0.6333 \\  
\hline
1/100&$d_{max}$  &0.0011  & 0.0015 &0.0015  &0.1667  &0.2500  & 0.2500 \\  
\hline                                                           
3/20&7 & 0.0066 & 0.0081 & 0.0030 & 0.3125 & 0.0005 & 0.3077 \\
3/20&$d_{max}$   &0.0158  & 0.0162 &0.0147  &0.4375  &0.7778  & 0.4737 \\  
\hline
3/30&7 & 0.0027 & 0.0025 & 0.0017 & 0.2727 & 0.4000 & 0.2857 \\
3/30&$d_{max}$   &0.0045  & 0.0058 &0.0040  &0.3462  &0.6250  & 0.4000 \\  
\hline
3/50&7 & 0.0007 & 0.0009 & 0.0008 & 0.2778 & 0.5454 & 0.2692 \\
3/50&$d_{max}$   &0.0026  & 0.0027 &0.0034  &0.3553  &0.7317  & 0.4694 \\  
\hline
3/100&$d_{max}$ &0.0014  & 0.0017 &0.0017  &0.1667  &0.2500  & 0.2500 \\  
\hline                                                           
6/20&$d_{max}$   &0.3102  & 0.1508 &0.3281  &0.4000  &0.6667  & 0.6000 \\  
\hline
6/30&$d_{max}$   &0.0576  & 0.0724 &0.0261  &0.4375  &0.7778  & 0.5882 \\  
\hline
6/50&$d_{max}$   &0.0079  & 0.0090 &0.0069  &0.4375  &0.7778  & 0.6471 \\  
\hline
6/100&7 & 0.0040 & 0.0047 & 0.0031 & 0.3125 & 0.4545 & 0.2857 \\
6/100&$d_{max}$  &0.0053  & 0.0042 &0.0050  &0.4130  &0.8125  & 0.5306 \\  
\hline                                                           
11/20&$d_{max}$  &0.0117  & 0.0137 &0.0131  &0.4000  &0.7500  & 0.5882 \\  
\hline
11/30&7 & 0.0033 & 0.0033 & 0.0029 & 0.2500 & 0.4000 & 0.3846 \\
11/30&$d_{max}$  &0.0052  & 0.0062 &0.0038  &0.3600  &0.6522  & 0.5455 \\  
\hline
11/50&7 & 0.0019 & 0.0031 & 0.0025 & 0.2308 & 0.5385 & 0.2778 \\
11/50&$d_{max}$  &0.0048  & 0.0055 &0.0042  &0.3684  &0.5500  & 0.4186 \\  
\hline
11/100&7 & 0.0009 & 0.0009 & 0.0014 & 0.2857 & 0.4667 & 0.4571 \\
11/100&$d_{max}$ &0.0034  & 0.0046 &0.0032  &0.3909  &0.8545  & 0.5397 \\  
\hline                                                           
14/20&$d_{max}$  &0.0251  & 0.0195 &0.0228  &0.4167  &0.7143  & 0.6364 \\  
\hline
14/30&7 & 0.0038 & 0.0057 & 0.0049 & 0.2500 & 0.4286 & 0.4000 \\
14/30&$d_{max}$  &0.0046  & 0.0056 &0.0041  &0.4250  &0.7391  & 0.5455 \\  
\hline
14/50&7 & 0.0037 & 0.0041 & 0.0051 & 0.2500 & 0.3333 & 0.4286 \\
14/50&$d_{max}$  &0.0067  & 0.0104 &0.0091  &0.3704  &0.6538  & 0.5909 \\  
\hline
14/100&7 & 0.0002 & 0.0007 & 0.0007 & 0.2000 & 0.6250 & 0.3226 \\
14/100&$d_{max}$ &0.0021  & 0.0028 &0.0022  &0.2576  &0.5373  & 0.3766 \\  
\hline
\end{tabular}\\
\vspace{0.2cm}
\label{compmetcap7dmax2}
\end{center}
\end{table}

As for the number of nondominated points by execution time ($\sigma^5$) and  by subproblem ($\sigma^6$) the results presented in Table \ref{compmetcap7dmax2} are similar to the one-dimensional case. The $FPA$ method obtained the greatest values in more than 75\% of the instances. However, in relation to $\sigma^5$ its variability  is similar to the other two methods and greater than the others considering $\sigma^6$ as shown in Figure \ref{boxplotPCE2D}.

Figure \ref{frontmet2Dc7} shows the Pareto front approximations obtained by the $FPA$, $LEC$, and $AWT$ methods considering $c=7$ for the  instance 11/100 of the $2D-BOCSP$. This instance was chosen because it has the greatest cardinality and hypervolume values (see Table \ref{sigma1-sigma4-2D-cap7dmax}). Note that the front obtained by the $FPA$ method, despite obtaining a smaller number of nondominated points, it has a greater hypervolume, and these points dominate almost all the points generated by the other two methods, thus providing solutions closer to the Pareto optimal  front. However, if there is interest in a greater number of points spread along the rectangle that contains the Pareto optimal front, the three methods should be combined as is further discussed in  Section \ref{Res3}.

 \subsection{Combining scalarization methods to obtain a Pareto Front Approximation}\label{Res3}
The computational results discussed in Sections \ref{Res1D} and \ref{Res2D}, draw attention to an important  question: ``What is the quality of the FrA obtained with the union of the points generated by the three scalarization methods (LEC, FPA and AWT)?". The answer to this question is explored in this section considering 51 out of the 60 instances presented in Sections \ref{Res1D} and \ref{Res2D} for which $\sigma_1 > 2$ for at least one of the three solution methods. 

The instances are divided into three groups. The  Group 1 is composed of 20 instances of the 1D-BOCSP, the Group 2  is composed of 13 instances of the 2D-BOCSP with $c=7$, and Group 3 is composed of 18 instances of the 2D-BOCSP with $c=dmax$.

\subsubsection*{Results for the instances of Group 1 (1D-BOCSP, $c=4,7, dmax$)}

For the Group 1, regarding cardinality, on average, the FPA method found 2\% more points than the LEC method. The AWT method found, on average, 31.1\% more than the LEC and FPA methods. As for the hypervolume metric, the FPA method performed, on average, 1.7\% better than the LEC method, and the LEC method performed, on average, 2.2\% better than the AWT method, while the FPA method performed, 4.7\% better than the AWT method. Therefore, regarding both metrics, cardinality and hypervolume, for the 1D case, the FPA method performed best in hypervolume and AWT in cardinality. In this sense, AWT generates more points, and FPA has better quality in the points found, as they are closer to the optimal boundary. 

For each instance, observing the relationship among the three methods, it is possible to see that there are points generated by one of the methods that are dominated by the points generated by the other two methods. However, no dominance pattern among the methods was identified. Considering the set of points obtained by the union of the points generated by the three methods, on average, 35.4\% of them would form an FrA, and the remaining 64.6\% were points equal to or dominated by the other points on the front. 

For each one of the 20 1D-BOCSP instances, the FrA formed by the union of the nondominated points was determined. It was noted that the hypervolume of this curve is, on average, 5.4\% greater than that found by the FrA generated by the LEC method, 7.2\% greater than that found by the FPA method, and 2.8\% greater than that found by the AWT method. Another important observation regarding cardinality is that this union of points forms a new FrA with an average cardinality 18\% greater than the ones generated by the LEC method, on average 40\% better than the FPA, and on average 22\% greater than the cardinality of the curve generated by the AWT method.

\subsubsection*{Results for the instances of Group 2 (2D-BOCSP, $c=7$)}

For the instances in Group 2, on average, the LEC method generated 15.8\% more points than the  FPA method and 3.2\% more than the AWT. Observing the relationship among the three methods, the AWT method found on average  13.8\% more points than the  FPA method. In relation to the hypervolume metric, the LEC method had a performance that was, on average, 7.0\% better than the FPA and 9.6\% better than the AWT method, while the AWT  method was 1.3\% better that the FPA on average. Then, taking into account both metrics, the LEC method was better than the two other methods.

As in the case for the 1D-BOCSP (instances in Group 1), Observing the relationship among the points generated for each instance, no dominance pattern was detected.  For the  set of points obtained by the union of the points generated by the three methods, on average, 35.0\% of them compose an FrA and the others (65.0\%) are the same or are dominated by the points of the new FrA.

For each one of the 13 instances a new FrA was constructed from the union of the points generated by the three methods. The hypervolume of this FrA is, on average, 4.0\% greater than the FrA generated by the LEC method, 15.2\% greater the FrA generated by the method FPA, and 14.0\% greater than the one generated by the AWT method. In relation to the cardinality, on average, these union of points generated an FrA 2.7\% greater than the one generated by the LEC method, 9.9\% greater than the one generated by the FPA, and  4.7\% greater than the one generated by  AWT method.

\subsubsection*{Results for the instances of Group 3 (2D-BOCSP, $c=dmax$)}

For the instances in Group 3, on average the LEC method generated 8.5\% more points than the  FPA method and 40.4\% more than the AWT. The FPA method found on average  30.2\% more points than the  AWT method. In relation to the hypervolume metric, the LEC method had, on average, a behaviour 3.6\% better than the FPA and the AWT methods, while the FPA method was 1.3\% better that the AWT on average. Then, taking into account both metrics, the LEC method was better than the two other methods.

As in the case for the 1D-BOCSP  and for  the 2D-BOCSP with $c=7$ (instances in  Groups 1 and  2 respectively), no dominance pattern was detected among the points found by the three different methods.  Regarding the  new set of points obtained by the union of the points generated by the three methods, 36.7\% of them compose an FrA and the others (63.3\%) are the same or are dominated by the points of the new FrA. 

For each one of the 18 instances an FrA was constructed from the union of the points generated by the three methods. The hypervolume of this new FrA is, on average, 0.6\% greater than the FrA generated by the LEC method, 1.0\% greater the the FrA generated by the method FPA, and 1.3\% greater than the one generated by the AWT method. In relation to the cardinality, on average, these union of points generated an FrA 1.0\% greater than the one generated by the LEC method, 6.1\% greater than the one generated by the FPA, and  37.8\% greater than the one generated by  AWT method.

Therefore, as in the other two Groups, it can be observed that the studied methods complement each other, as the union of the different generated points significantly improves the cardinality and hypervolume of the FrA. In this context, an alternative solution is to  make a union of the points generated by the three methods and take the nondominated points from this set to compose a new FrA.

\section {Final Remarks}\label{concl}

The aim of this paper is to propose and solve a mathematical optimization model for the bi-objective cutting stock problem. The two objectives considered are the minimization of the total number of objects and the minimization of  the total number of saw cycles. Three scalarizations methods from the literature are discussed and adapted  to solve two versions of the problem, the one-dimensional and the two-dimensional cases. A dynamic  column generation approach is employed. 

In the literature review presented, it was shown that this problem has been scarcely studied, most articles address the problem disconsidering its multi-objective nature.  A summary of the main concepts of multi-objective optimization and a complete description of the mathematical optimization model proposed is presented, including the procedures employed to generate columns for the one-dimensional and the two-dimensional versions of the problem. The reasons for choosing the classical Lexicographic Constraint method (LEC),  Front Partitioner Algorithm (FPA) and the  the Augmented Weighted Tchebycheff Method (AWT) are explained and details of the adaptations necessary to employ them to solve the BOCSP with column generation are also given.

In the computational experiment, at first it is shown the advantages of employing a dynamic column generation approach. This study was conducted for the 1D-BOCSP, with the saw capacity equals to 7. Considering the all three methods together, the result shows that the Dynamic column generation provides an FrA closer to the Pareto optimal front than the Static Column Generation, with an average of  10.81\% increase in the total computational time. Therefore the Dynamic column generation was employed in the remainder of the computational study.

The performances of the three scalarization methods (LEC, FPA and AWT) implemented with Dynamic column generation to solve the BOCSP were analysed. The trade-off between the minimization of the total number of objects and minimization of the total number of saw cycles was studied using 51 instances with different characteristics related to the problem dimension, type and number of items. Considering the metrics cardinality and hypervolume, the LEC method presented, on average, the best performance for the two-dimensional instances, followed by the FPA and AWT methods. For the one-dimensional instances the FAP method had a better performance related to hypervolume metric and the method AWT had a better performance related to the cardinality metric.

It was observed  that some of the points generated by one method dominate the points generated by another method. However, no pattern of dominance was detected. On the other hand, it was determined that when joining all the points generated by the three methods, on average, 36\% of them would compose a new FrA, while  the remaining 64\% are either identical (repeated) points or dominated by the points of the new curve. The cardinality and hypervolume of this new FrA are always greater than those generated by the methods studied. In this sense, the new curve has better quality than the others generated by the previous methods. In this context, the LEC, FPA and AWT methods can be used to complement each other.

\section*{Funding}
Research partially supported by the National Council for the Improvement of Higher Education  - CAPES  (Coordenação de Aperfeiçoamento de Pessoal de Nível Superior); by the São Paulo Research Foundation - FAPESP (grants  2022/05803-3 and  2013/07375-0); and by the National Council for Scientific and Technological Development  - CNPq  (Conselho Nacional de Desenvolvimento Científico e Tecnológico){(grant 306518/2022-8)}.
\clearpage
\bibliographystyle{abbrv}

\bibliography{B-referencias}

\end{document}